\documentclass[12pt,a4paper,final]{article}% scrartcl}

\usepackage[margin=24.4mm]{geometry}
\usepackage[utf8]{inputenc}  
\usepackage{xcolor,scrtime,graphicx}
\usepackage{amsmath,amsfonts,amssymb,mathrsfs}

\usepackage{bef_alex,bm} %%PWD (line used to be absent)
\def\FG{\bm}
\def\mafo{\mathrm}

\numberwithin{equation}{section}
\numberwithin{figure}{section}
\usepackage{todonotes} 
 
 \newcommand{\TTODO}[1]{}

\usepackage[notref,notcite]{showkeys}

\usepackage{tikz}

\makeatletter
\renewcommand*\env@cases[1][1.2]{%
  \let\@ifnextchar\new@ifnextchar
  \left\lbrace
  \def\arraystretch{#1}%
  \array{@{\,}c@{\quad}l@{}}%
}
\makeatother

\newtheorem{problem}[theorem]{Problem}

\newcommand{\mfD}{\mathfrak D}

\newcommand{\mfP}{\mathfrak P}
\newcommand{\sfw}{\mathsf w}
\newcommand{\sfE}{\mathsf E}
\newcommand{\sfR}{\mathsf R}
\newcommand{\sfC}{\mathsf C}

\renewcommand{\ti}{{\times}}
\newcommand{\scrK}{\mathscr K }
\newcommand{\scrL}{\mathscr L}
\newcommand{\scrB}{{\mathscr B}}
\newcommand{\frakB}{{\mathfrak B}}

\newcommand{\weak}{\rightharpoonup}

\newcommand{\Gto}{\overset{\Gamma}\to}

\newcommand{\dom}{\mathop{\mathrm{dom}}} 
\newcommand{\eff}{\mathrm{eff}} 
\newcommand{\red}{\mathrm{red}} 
\newcommand{\fast}{\mathrm{fast}} 
\newcommand{\slow}{\mathrm{slow}} 
\newcommand{\LB}{\lambda_\rmB} 

\newcommand{\SI}{\mafo{SI}} 
\newcommand{\IS}{\mafo{IS}} 

\newcommand{\PPO}{\bbP}
\newcommand{\PPC}{\bbP^*}

\newcommand{\STEP}[1]{\noindent{\itshape #1}} 

\begin{document}
 
\title{Non-equilibrium steady states as 
saddle points and\\
EDP-convergence for slow-fast gradient systems}

\author{Alexander Mielke\\ 
WIAS Berlin and  Humboldt-Universit\"at zu Berlin}

\date{9 March 2023, revised 10 October 2023} 

\maketitle

{\def\thefootnote{}
\footnotetext{e-mail: \texttt{alexander.mielke@wias-berlin.de}}}

\begin{abstract} The theory of slow-fast gradient systems leads in a natural
  way to non-equilibrium steady states, because on the slow time scale the fast
  subsystem stays in steady states that are controlled by the interaction with
  the slow system. Using the theory of convergence of gradient systems
  depending on a small parameter $\eps$ (here the ratio between the slow and
  the fast time scale) in the sense of the \emph{energy-dissipation principle}
  shows that there is a natural characterization of these non-equilibrium
  steady states as saddle points of a so-called $B$-function where the slow
  variables are fixed.  We give applications to slow-fast reaction-diffusion
  systems based on the so-called cosh-type gradient structure for reactions. It
  is shown that two binary reaction give rise to a ternary reaction with a
  state-dependent reaction coefficient. Moreover, we show that a
  reaction-diffusion equation with a thin membrane-like layer convergences to a
  transmission condition, where the formerly quadratic dissipation potential
  for diffusion convergences to a cosh-type dissipation potential for the
  transmission in the membrane limit.
\end{abstract}

%{\def\thefootnote{}
%\footnotetext{arXiv:2303.07175 (v1 13.3.2023) \hfill Datum \today, Uhrzeit \thistime\ h}}

%{%\def\tocname{}
%\footnotesize\tableofcontents}

\section{Introduction}
\label{se:Intro}

\subsection{Dissipative evolution equations and gradient systems}
\label{suu:DissGS}

We consider systems that are characterized by a state $u$ in the state space
$X$ (a convex subset of a Banach space or a manifold).  The evolution process
of the system is assumed to be slow enough such that inertial effects can be
neglected. Moreover, we assume that the dynamics can be described by a balance
of friction forces $\xi$ and a potential restoring forces $-\rmD\calE(u)$,
where the energy potential $\calE$ can be a free energy or a negative
entropy. The friction forces are given by a kinetic relation
\[
\xi = K(u,\dot u),  
\]
where the mapping $K(u,\cdot) : X\to X^*$ is assumed to satisfy $K(u,0)=0$ and
$\langle K(u,v),v\rangle \geq 0$ for thermodynamic reasons. Thus, the evolution
equations of our interest are 
given by 
\[
0= K(u,\dot u) + \rmD\calE(u). 
\]  
Classical viscous friction leads to a linear kinetic relation $\xi = K(u,\dot
u) = \bbK(u) \dot u$, while Coulomb friction satisfies $\xi = \mu(u)\,\mafo{Sign}(\dot
u)$. Along solutions we obtain decay of the (free) energy via
\[
\frac\rmd{\rmd t} \calE(u(t)) = \langle \rmD\calE(u(t)), \dot u(t)\rangle = -
\langle  K(u(t),\dot u(t)), \dot u(t)\rangle \leq 0. 
\]

For general nonlinear kinetic relations, the proper generalization of Onsager's
symmetry $\bbK(u)=\bbK(u)^*\geq 0$ for linear kinetic relations (see
\cite{Onsa31RRIP1}) is the existence of a so-called \emph{dissipation
  potential} $\calR$ such that $K(u,v)=\rmD_v\calR(u,v)$, see the discussion in
\cite{MiRePe16GORR}.  A function $\calR: X\ti X \to [0,\infty]$ is called a
dissipation potential, if for all $u\in X$, the function
$\calR(u,\cdot): X \to [0,\infty]$ is lower semi-continuous, convex and
satisfies $\calR(u,0)=0$. The dual dissipation potential
$\calR^*: X\ti X^* \to [0,\infty]$ is defined via
\[
\calR^*(u,\xi) :=\sup\bigset{\langle \xi,v\rangle - \calR(u,v)}{ v\in X }.
\] 
By classical duality theory (the so-called Fenchel equivalences
\cite{Fenc49CCF}) we have 
\[
\xi = \rmD_v \calR(u,v) \quad \Longleftrightarrow \quad v=\rmD_\xi \calR^*(u,\xi)  
 \quad \Longleftrightarrow \quad 
 \calR(u,v)+\calR^*(u,\xi) \leq  \langle \xi,v\rangle.
\]

The triple $(X,\calE,\calR)$ will shortly be called a gradient system (GS), and
the gradient-flow equation associated with the GS $(X,\calE,\calR)$ can be
written in two equivalent forms, namely
\begin{equation}
  \label{eq:I.GFE}
  \text{(I) \ } 0 = \rmD_v\calR(u,\dot u)+ \rmD\calE(u) \quad \Longleftrightarrow \quad 
\text{(II) \ } \dot u = \rmD_\xi\calR^*\big(u,{-}\rmD\calE(u)\big). 
\end{equation}
Throughout this work we assume
for simplicity that $\calR$ and $\calR^*$ are differentiable; however most of
the theory carries over to the general case of convex functions where
$\rmD_v\calR$ and $\rmD_\xi\calR^*$ are replaced by the set-valued convex
subdifferentials, see e.g.\ \cite{Miel23IAGS}. 

The restriction to evolution processes given by GS restricts the dynamics to
dissipative processes only, but leads to several favorable thermodynamical
principles. A first observation is that all steady states $u_*$ of
\eqref{eq:I.GFE} are equilibria, i.e.\ $\rmD \calE(u_*)=0$. One purpose
of this work is to show that in slow-fast GS one can have non-equilibrium
phenomena such as Non-Equilibrium Steady States (NESS) if the slowly changing
components of the system allow the very fast changing components to relax into
a steady state but not into an equilibrium. 

Next, we will show how gradient-flow equations can be reformulated in a
thermodynamic variational form, and then explain how the limit of slow-fast
systems can be derived. On this way, NESS will appear naturally and a new
saddle-point  characterization of NESS will emanate. 

\subsection{Energy-dissipation principle for slow-fast systems}
\label{su:EDP.slowfast}

Under suitable technical assumptions, the gradient-flow equation
\eqref{eq:I.GFE} can be reformulated in a thermodynamical way by the so-called
\emph{energy-dissipation principle} (EDP), see
\cite{RoMiSa08MACD, Miel23IAGS}.  If $u:[0,T]\to X$ satisfies
$\int_0^T \big[\calR(u,\dot u){+}\calR^*(u,{-}\rmD \calE(u)\big) \big] \dd t <
\infty$, then (I) and (II) hold a.e.\ in $[0,T]$ if and only if the
\emph{energy-dissipation inequality} (EDI) holds, namely
\begin{equation}
  \label{eq:I.EDI}
  \calE(u(T)) + \int_0^T\!\! \Big(\calR(u,\dot u){+}\calR^* \big(u,{-}\rmD
\calE(u)\big) \Big) \dd t \leq \calE(u(0)). 
\end{equation}

If we now have a family $(X,\calE_\eps,\calR_\eps)_{\eps>0}$ of GS with a small
parameter $\eps>0$, we say that this family converges in the sense of the EDP
to the limiting GS $(X,\calE_\eff,\calR_\eff)$ if we have the following
$\Gamma$-convergences 
\begin{equation}
  \label{eq:I.def.EDPcvg}
  \calE_\eps \Gto \calE_\eff \text{ in } X \quad \text{and} \quad  \mfD_\eps
  \Gto \mfD_0 \text{ \ in } \rmL^2([0,T];X),
\end{equation}
where the dissipation functionals $\mfD_\eps:\rmL^2([0,T];X\to [0,\infty]$ are
defined as follows:
\begin{align*}
&\mfD_\eps\big(u(\cdot)\big) := \int_0^T\!\! \Big(\calR_\eps(u,\dot
u){+}\calR^*_\eps \big(u,{-}\rmD \calE_\eps(u)\big) \Big) \dd t
&&\text{for }\eps>0,
\\ 
&\mfD_0\big(u(\cdot)\big) := \int_0^T\!\! \Big(\calR_\eff(u,\dot
u){+}\calR^*_\eff \big(u,{-}\rmD \calE_\eff(u)\big) \Big) \dd t 
&&\text{for }\eps=0. 
\end{align*}
We refer to \cite{LMPR17MOGG} for the first discussion of this concept, to
\cite{MiMoPe21EFED} for refinements, and to \cite{DoFrMi19GSWE, Fren19DEGS,
  MiPeSt21EDPC, FreLie21EDTS, PelSch22?CGST} for various applications of this
approach. 

We emphasize two important properties of EDP-convergence: The first one simply
states that if $u_\eps$ are solutions to $(X,\calE_\eps,\calR_\eps)$, we have
convergence $u_\eps (t) \to u(t)$ in a suitable way, and the initial conditions
$u_\eps(0)$ are well-prepared, i.e.\
$\calE_\eps(u_\eps(0)) \to \calE_\eff(u(0))$, then $u$ is a solution of the
effective GS $(X,\calE_\eff,\calR_\eff)$. The second property states that
$\calR_\eff$ can be different to a potentially existing $\Gamma$-limit
$\calR_0$, i.e.\ $\calR_\eps \Gto \calR_0$. The point is that $\mfD_\eps$
involves a nonlinear construction for the pair $(\calE_\eps,\calR_\eps)$, which
allows for a transfer of microscopic information of the energy (encoded in
$\calE_\eps{-}\calE_\eff$) into the dissipation $\calR_\eff$. We will see this
below in Section \ref{su:DiffMembrOtt} where $\calR_\eps^*(u,\cdot)$ is
quadratic and has a quadratic limit $\calR^*_0$ but $\calR_\eff^*$ contains a
cosh-type membrane part for the transmission through the membrane.

The slow-fast GS under consideration are assumed to be of the following form 
\[
  X=X_\slow\ti X_\fast, \quad \calE_\eps(U,w)=E(U)+ \eps \,e(w), \quad
  \calR^*_\eps(U,w;\Xi, \xi) = \ol\calR^*\!\big(U,w;\Xi,
  \tfrac{\ds1}{\ds\eps}\;\!\xi\big),
\]
which provides only one class of GS where slow-fast effects can be studied (see
\cite{MieSte20CGED, MiPeSt21EDPC} for other scalings). The associated
gradient-flow equation reads
\begin{equation}
  \label{eq:I.FastSlowSy}
  \binom{\ds \dot U}{\ds \eps\;\!\dot w} = \pl_{\Xi,\zeta} 
\ol\calR^*\!\big(U,w; -\rmD E(U), - \rmD e(w)\big),
\end{equation}   
which shows nicely the slow-fast structure, because $\eps$ only appears once,
namely in front of the time derivative $\dot w$ of the fast variable
$w \in X_\fast$. 

As in \cite{Kueh15MTSD} the above two-scale system gives rise
to the ``slow system'' (by setting formally $\eps=0$) 
\begin{equation}
  \label{eq:I.SlowSyst}
   \binom{\ds \dot U}{0} = \pl_{\Xi,\zeta} 
\ol\calR^*\!\big(U,w; -\rmD E(U), - \rmD e(w)\big),
\end{equation}
and to the ``fast system'' (on the short time scale $s = t/\eps$ with $'=
\rmd/\rmd s = \eps \rmd /\rmd t$) 
\begin{equation}
  \label{eq:I.FastSyst}
   U=\text{constant}, \quad w' = \rmD_\zeta \ol\calR^*\!\big(U,w; -\rmD E(U), -
   \rmD e(w)\big). 
\end{equation}
The steady states $w=\wt\sfw(U)$ of \eqref{eq:I.FastSyst} can be seen as NESS
that are nontrivial because of the given fixed values of $U$. In Section
\ref{se:ConstrSP} we discuss our first main goal, namely a variational,
thermodynamical characterization of NESS as saddle points. The second main
goal is to show that the reduced slow system, obtain by inserting
$w=\wt\sfw(U)$ into \eqref{eq:I.SlowSyst} is again a gradient-flow equation for
a GS $(X_\slow, \calE_\eff, \calR_\eff)$ for $\calE_\eff =E$ and a suitable
effective dissipation potential $\calR_\eff$, see Section
\ref{se:ReductGSNESS}.

As a simple example we may consider the quadratic GS $(X_\slow\ti X_\fast,
\calE_\eps, \calR_\eps)$ with 
\[
\calE_\eps(U,w)=\frac12\langle \bbA_\rms U,U\rangle + \frac\eps2 \langle
\bbA_\rmf  w,w\rangle \ \text{ and } \ 
\calR^*_\eps(\Xi,\zeta)= \frac12 \langle \bbK_{\rms\rms}\Xi,\Xi\rangle
+ \frac1\eps \langle \bbK_{\rms\rmf}\zeta,\Xi\rangle + \frac1{2\eps^2} \langle
\bbK_{\rmf\rmf} \zeta,\zeta\rangle.   
\]
The associated slow system reads 
\[
\binom{\dot U}0 = \binom{\bbK_{\rms\rms} \bbA_\rms U + \bbK_{\rms\rmf} \bbA_\rmf w} 
{\bbK_{\rms\rmf}^* \bbA_\rms U + \bbK_{\rmf\rmf} \bbA_\rmf w}  \quad
\Longleftrightarrow \quad  
\begin{cases} \dot U= \bbK_\eff \bbA_\rms U \text{ with } \bbK_\eff= \bbK_{\rms\rms}{-}
  \bbK_{\rms\rmf} \bbK_{\rmf\rmf}^{-1} \bbK_{\rms\rmf}^*\\
\text{and } w = - (\bbK_{\rmf\rmf}\bbA_\rmf)^{-1} \bbK_{\rms\rmf}^* \bbA^{}_\rms U.
\end{cases} 
\]
For this example, it is easy to find the NESS
$ w = \wt\sfw(U)= - (\bbK_{\rmf\rmf}\bbA_\rmf)^{-1} \bbK_{\rms\rmf}^*
\bbA^{}_\rms U$ as well as the effective GS $(X_\slow, \calE_\eff,\calR_\eff)$ with
$\calE_\eff(U)= \frac12 \langle \bbA_\rms U,U\rangle$ and
$\calR^*_\eff(\Xi)=\frac12 \langle \Xi, \bbK_\eff \Xi\rangle$, see Section
\ref{su:QuadrEnerDiss} for more details. The purpose of
this work is to generalize these trivial observations to more general settings with
\emph{state-dependent} and \emph{nonlinear} kinetic relations.

\subsection{Reduction via B-functions and NESS}
\label{su:B-function.SaPo}

To see where the motivation for out theory comes from, we return to the
two-scale equation \eqref{eq:I.FastSlowSy} and apply the EDP convergence theory
as in \eqref{eq:I.def.EDPcvg}. We observe that $\mfD_\eps$ takes the
simple form
\[
\mfD_\eps\big(U,w\big) := \int_0^T\!\! \Big(\ol\calR(U,w;\dot U, \eps\!\;\dot
w) {+}\ol\calR^*\!\big(U,w;{-}\rmD E(U), -\rmD e(w)\big) \Big) \dd t, 
\]
and it is tempting to drop the term $\eps\dot w$ and minimize the integrand
for each $t\in [0,T]$ with respect to $w$. However, we will see that this
approach is not correct because we have to find the correct
\emph{non-equilibrium steady states} which create a nontrivial flux as a limit
of $\eps \dot w$. 

We follow the approach in \cite{LMPR17MOGG} and estimate $\ol\calR$ from below
via
\[
\ol\calR(U,w;\dot U, \eps\!\;\dot
w) \geq \Big\langle \binom{\Xi}{\zeta},\binom{\dot U}{\eps\!\;\dot
w} \Big\rangle - \ol\calR^*(U,w;\Xi,\zeta),
\]
where $(\Xi,\zeta):[0,T]\to X_\slow^* \ti X_\fast^*$ are smooth test
functions. Thus, we have 
\begin{align*}
&  \mfD_\eps(U,w)\geq \int_0^T \Big(\Big\langle \binom{\Xi}{\zeta},\binom{\dot
    U}{\eps\!\;\dot w} \Big\rangle - \frakB_{\ol\calE,\ol\calR}(U,w;\Xi,\zeta)
  \Big) \dd t,
\end{align*}
where $\ol\calE(U,w)=E(U)+e(w)$. For general GS $(X,\calE,\calR)$ we
call 
\[
\frakB_{\calE,\calR}: \mafo{dom}(\frakB_{\calE,\calR}) \subset X\ti X^* \to \R;\
(u,\xi)\mapsto \calR^*(u,\xi) - \calR^* \big(u,{-}\rmD\calE(u)\big) 
\]
the associated ``B-function''.  

It turns our that we can now pass to the limit $\eps\to 0$ by omitting the term
$\eps\!\;\dot w$, because it multiplies the given smooth test function
$\zeta$. Then, we can maximize with respect to $\zeta$ and
minimize with respect to $w$ for each individual $t \in
[0,T]$. Hence, in terms of the 
function $\frakB_{\ol\calE,\ol\calR}$ we are lead to the following sup-inf
problem:
\begin{equation}
  \label{eq:I.Lred}
  \scrB_\red(U,\Xi):= \sup_{w\in X_\fast} \inf_{\zeta\in X^*_\fast}
\frakB_{\ol\calE,\ol\calR}(U,w;\Xi,\zeta) .
\end{equation}
We could have equally well defined $\scrB_\red$ in terms of the
corresponding inf-sup problem, which in general would produce a larger
function, see \eqref{eq:EstimSP.a}. However, our theory is based on the
existence of global saddle points, which imply that the sup-inf equals the
inf-sup. Our choice of taking the infimum first is motivated by the simplicity
of minimizing the convex functional
$\zeta \mapsto \frakB_{\ol\calE,\ol\calR}(U,w;\Xi,\zeta)$ over
$\zeta \in X^*_\fast$.  By construction, the function
$(\Xi,\zeta) \mapsto \frakB_{\ol\calE,\ol\calR}(U,w;\Xi,\zeta)$ is always
convex. The classical existence theory for global saddle points works under the
additional assumption that
$(U,w) \mapsto \frakB_{\ol\calE,\ol\calR}(U,w;\Xi,\zeta)$ is concave.  This
condition is valid for Otto's gradient structure for diffusion (cf.\ Section
\ref{su:DiffMembrOtt}), but not for nonlinear chemical reaction systems as
treated in Section \ref{suBinaryTernaryReact}.
For our general theory, we do not need to assume concavity of 
$(U,w) \mapsto \frakB_{\ol\calE,\ol\calR}(U,w;\Xi,\zeta)$.

We say that the reduced B-function $\scrB_\red$ has \emph{BER
  structure} if there exists an effective dissipation potential
$\calR_\eff : X_\slow\ti X_\slow \to [0,\infty]$ such that it can be written as
\begin{equation}
  \label{eq:I.DualStruct}
  \scrB_\red (U,\Xi) = \frakB_{E,\calR_\eff}(U,\Xi) 
  = \calR_\eff(U,\Xi) - \calR^*_\eff\big(U,{-}\rmD E(U)\big).
\end{equation}
Using the EDP backwards, we see that the effective GS
$(X_\slow,E,\calR_\eff)$ with the gradient-flow equation 
\[
 \dot U  = \pl_\Xi \calR^*_\eff \big(U, -\rmD \calE (U)\big)
\]
indeed describes the limiting dynamics. 

Thus, the main point in applying this theory successfully is to show the
existence of the BER structure $(E,\calR_\eff)$ for the reduced B-function
$\scrB_\red$. And it is here were the theory of NESS comes into play.  The
definition of NESS in the above context means that we fix $\ol U\in X_\slow$ and
want to find the NESS $\ol w\in X_\fast$ such that 
\begin{equation}
  \label{eq:I.NESS}
  \binom00 = 
\binom{\rmD_\Xi \ol\calR^*\!\big( U,\ol w; -\rmD E(U), - \rmD  e(\ol w)\big)}
{\rmD_\zeta \ol\calR^*\!\big( U,\ol w; -\rmD E(U), - \rmD  e(\ol w)\big)}
 - \binom{\ol   y}{0} , \quad U=\ol U \in X_\slow, \quad \ol y\in X_\slow.
\end{equation}
We refer to \eqref{eq:NESSeqn} for the general case involving port mappings
$P:X\to Y$ and $\PPC:X^*\to Y^*$, which reduce to \eqref{eq:I.NESS} if we
choose $Y=X_\slow$, $P(U,w)=U$ and $\PPC(\Xi,\zeta)=\Xi$. We observe that
fixing $U=\ol U$ artificially generates a flux $\ol y$ which is generated by
the NESS $\ol w$ associated with $\ol U$, i.e.\ one can think of $\ol y$ as a
Lagrange multiplier for the constraint $U=\ol U$.

\subsection{Global saddle points provide NESS}
\label{su:GLoSAPoNESS}
 
The first major link between the theory of NESS and the above saddle-point
reduction for B-functions is the fact that a NESS $\ol w$ solving
\eqref{eq:I.NESS} gives rise to a \emph{global null-saddle}
$(w,\zeta) = \big(\ol w, -\rmD e(\ol w)\big)$ for
$\frakB_{\ol\calE,\ol\calR}\big(\ol U, \cdot \ \! ;-\rmD E(\ol U),\cdot\ \!
\big)$, i.e.\ for all $(w,\zeta) \in X_\fast\ti X_\fast^*$ we have 
\begin{equation}
  \label{eq:I.GlobSaddle}
  \begin{aligned}
& 
\frakB_{\ol\calE,\ol\calR}\big(\ol U, w;-\rmD E(\ol U),-\rmD e(\ol w) \big)
\leq 0 
\\
& \qquad =
\frakB_{\ol\calE,\ol\calR}\big(\ol U, \ol w;-\rmD E(\ol U),-\rmD e(\ol w) \big)
 \leq \frakB_{\ol\calE,\ol\calR}\big(\ol U, \ol w;-\rmD E(\ol U),\zeta \big).
\end{aligned}
\end{equation}
Section \ref{su:ConstrainedSP} discusses conditions under which null-saddles
$\big(\ol w,\zeta)$ automatically satisfy $\zeta=-\rmD e(\ol w)\big)$, where $\ol w $ is
NESS solving \eqref{eq:I.NESS}. There seem to exist a number of different
variational characterizations (also called extremum principles) of NESS, but to
the best of the author's knowledge the saddle-point formulation given here is
new. We refer to \cite{SaGaBr95VVPF, StrWei98MLEP, DDMU12RNSR},
\cite[Cha.\,V]{DegMaz84NET}, and \cite[Ch.\,30, pp.\,213-215]{Tsch00FESS}. In
particular, \cite{StrWei98MLEP} has the appealing title ``{\it Maximum of the
  Local Entropy Production Becomes Minimal in Stationary Processes}''.

The second important link arises from the fact that the existence of
null-saddles implies $\scrB_\red \big(\ol U,{-}\rmD \calE(\ol U) \big)=0$ with
$\scrB_\red$ from \eqref{eq:I.Lred}.  However, Proposition
\ref{pr:CondsDualStruc} shows that this condition (for all $\ol U\in X_\slow$)
is exactly the crucial condition for the existence of a BER structure in
the sense of \eqref{eq:I.DualStruct}. Theorem \ref{th:Lred.DualityStr}
provides the main result giving the explicit construction of $\calR_\eff$
in the form $\calR_\eff(U,\Xi)=\scrB_\red(U,\Xi)-\scrB_\red(U,0)$.

Section \ref{se:ReductGSNESS} gives a more detailed account of the reduction of
slow-fast GSs as discussed above. In particular, Section
\ref{su:Case2Product} also treats the case where the slow component
$U\in X_\slow$ and the fast component $w\in X_\fast$ only interact by a
constraint $P_\slow U=P_\fast w$, where $P_\slow:X_\slow \to Y$ and
$P_\fast : X_\fast \to Y$ are linear port mappings. In that case
the effective dual dissipation potential is the sum
\[
\calR^*_\eff(U, \Xi) = \calR^*_\slow (U, \Xi)  + \sfR_Y^*\big( P_\slow U,
\PPC_\slow \Xi\big),
\]
i.e.\ $\sfR_Y$ encodes all necessary information on the NESS in $X_\fast$.

\subsection{Applications in reactions and diffusion}
\label{su:ApplReacDiff}

Section \ref{se:ODEExamples} provides two ODE examples, the first being that of
a general quadratic dissipation potential and quadratic energies $E$ and
$e$. Everything can be explicitly calculated such that this case is helpful to
obtain guidance when the abstract theory may be overwhelming. The second
example treat a reaction-rate equation for four species $A$, $B$, $C$, and $D$
undergoing two binary reaction pairs $ A+ B \rightleftharpoons D$ and
$ A+ D \rightleftharpoons C$. Starting with constant reaction coefficients
$\kappa_{1,2}$ for the two reactions and assuming that the vector of equilibrium
densities is $(a_*,b_*,c_*,d_\eps)$ with $d_\eps = \eps w_\eps$ the
transformation $d(t)=\eps w(t)$ provides exactly a slow-fast GS as above, where
the energies $E$ and $e$ are relative Boltzmann entropies and $\ol\calR^*$ is
of cosh-type. Applying the above reduction method via NESS we find an effective
GS of cosh-type for the density vector $(a,b,c)$ that corresponds to the
ternary reaction pair $ 2A+ B \rightleftharpoons C$. The interesting point is
that, in contrast to the result in \cite{MiPeSt21EDPC}, the cosh-type gradient
structure is preserved, but now the effective reaction coefficient depends on
the density $a$.

Section \ref{se:LinearDiffMembrane} revisits the results obtained in
\cite{LMPR17MOGG} but now from a more general perspective. Moreover, the
results are generalized by allowing for a reaction term which models sorption
into and desorption from the background. The model starts from a
one-dimensional diffusion on an interval, where the diffusion coefficient in
the central membrane region ${]{-}\eps,\eps[}$ is scaled by $\eps$. Using Otto's gradient
structure (see \cite{Otto96DDDE, Otto98DLPF, JoKiOt98VFFP, Otto01GDEE}) we
start again from relative Boltzmann entropies $E$ and $e$ and from quadratic 
dual dissipation potentials $\ol\calR^*(U,w;\,\cdot\,,\,\cdot\,)$. In the limit
$\eps \to 0$ the membrane part collapses to an interface generating
transmission conditions. Our methods shows that $\sfR_Y$ is of cosh-type, which
indicates that it has inherited properties from the Boltzmann entropy $e$. Indeed,
the Boltzmann function $\LB(z)= z\log z -z +1$ with $\LB'(z)=\log z$ generates
by the saddle-point problem the cosh-type function $\sfC^*(\zeta)=
4\cosh(\zeta/2)-4$. Theorem \ref{th:MembrFormula} contains a much shorter
derivation of $\calR_\eff$ than in \cite{LMPR17MOGG, PelSch22?CGST}, and
Theorem \ref{th:MembRDS} generalizes the result to the case including a
reaction term that scales like $1/\eps$ in the membrane region
${]{-}\eps,\eps[}$. 

Finally, Appendix \ref{se:Saddles} provides the classical result on global saddle
points as discussed in \cite{EkeTem74ACPV}. For the readers convenience, we
include a full proof for the existence of saddle points for convex-concave
B-functions.

\section{NESS and constrained saddle points}
\label{se:ConstrSP}

We first provide our definition of NESS for general \emph{port gradient
  systems} and derive a few fundamental properties. Next we collect some
basic facts about unconstrained saddle points, then introduce the notion of
constrained saddle points using linear port mappings $P:X\to Y$ and
$\PPC : X^*\to Y^*$ and show that under a suitable additional condition that
these constrained saddle points are indeed NESS. Section
\ref{su:NESS.minimizer} shows a further characterization of NESS as
null-minimizers of a suitable auxiliary functional. Section
\ref{su:ReducLagrDual} provides the main result concerning the BER structure
for reduced B-functions $\scrB_\red$ if the associated NESS are null-saddles.

\subsection{Port gradient systems and NESS}
\label{su:PGS.NESS}

For the definition of NESS, we generalize the setting from slow-fast GS to
general port gradient systems (PGS) $(X,\calE,\calR,\PPO)$, where the linear
port mapping $\PPO:Y\to X$ maps the port space $Y$ to the (tangent space of)
$X$. The dual mapping $\PPC:X^*\to Y^*$ maps the thermodynamical driving forces
in $X^*$ to the reduced ones in $Y^*$. In the slow-fast GS we choose
$Y=X_\slow$ and $\PPO y=(y,0) \in X=X_\slow \ti X_\fast$. The ``fast system''
\eqref{eq:I.FastSyst} generalizes to the port gradient-flow equation
\begin{equation}
  \label{eq:PGS.GFE}
  \dot u= \rmD_\xi\calR^*\big(u, - \rmD\calE(u)\big) - \PPO y \quad  \text{with }
  y\in Y \ \text{ and } \ \PPC \rmD\calE(u)=-\eta \in Y^*. 
\end{equation}
Like in the theory of port-Hamiltonian systems (cf.\
\cite{EbMaVa07EHST,VanMeh23LPHD}), one can consider $y(t) \in Y$ as an input
(possibly time-dependent) and $\eta(t) \in Y^*$ as the output, or vice versa. A
simple calculation gives, along solutions, the power balance
\[
\frac\rmd{\rmd t} \calE(u) = -\underbrace{\big\langle {-} \rmD\calE(u),
  \calR^*(u,{-}\rmD \calE(u)) \big\rangle_{X}}_{\geq 0} +\langle \eta,y\rangle_Y,
\]
see e.g.\ \cite[Eqn.\,18]{EbMaVa07EHST}. 

\begin{definition}[NESS for a PGS] 
A state $\ol u\in X$ is called a NESS for the PGS $(X,\calE,\calR,\PPO)$ if
there exist a constant pair $(y,\eta)\in Y\ti Y^*$ such that $u(t)=\ol u$ is a solution of
\eqref{eq:PGS.GFE} and $\langle \eta, y\rangle_Y \neq 0$. This means that $\ol
u$ is a NESS if and only if it satisfies the NESS equations
\begin{equation}
  \label{eq:NESSeqn}
   0=\rmD_\xi\calR^*(\ol u,-\rmD\calE(\ol u)) - \PPO  y, \qquad
    \PPC\rmD\calE( \ol u)= - \eta \in Y^*, \qquad \ol y \in Y. 
\end{equation} 
\end{definition} 
 
Defining the dissipation function $\Phi_*(u,\xi):= \langle \xi, \rmD_\xi
\calR^*(u,\xi)\rangle_X \geq 0$ we see that every NESS dissipates energy 
\[
\langle \eta,y\rangle_Y = \Phi_*(\ol u, {-}\rmD\calE(\ol u))>0.
\]
This is in contrast to equilibrium solutions $\wt u\in X$ satisfying
$0=\rmD_\xi \calR^*(\wt u,{-}\rmD\calE(\wt u))$.\medskip

Under the assumption that for all $\eta\in Y^*$ there exists a unique NESS
$\ol u_\eta$ of \eqref{eq:NESSeqn} with Lagrange parameter $y=\ol y_\eta$, we
can define the \emph{port relation}
\[
\mfP:Y^* \to Y; \ \eta \mapsto \ol y_\eta.
\]
Such port relations which will play a crucial role in the sequel. As a
first result we observe, that in the case that $\calR$ is independent of the
state, the port relation can be obtained easily from $\calR$, it is independent
of the energy $\calE$, and it is given as the differential of an effective
potential $\sfR_Y$. We refer to Section \ref{su:QuadrEnerDiss} for a simple
and explicit case.

\begin{proposition}[Port relation for state-independent dissipation]
\label{pr:PortRelaSIDiss}
If $\calR:X\to [0,\infty]$ is a state-independent dissipation potential, then
the  port relation $\mfP$ is given by 
\[
y=\mfP(\eta)= \pl \sfR^*_Y (\eta) \quad \text{with } \ 
\sfR_Y(v)=\calR(\PPO v).
\]
Equivalently, $\sfR^*_Y$ is characterized via $\sfR^*_Y(\eta):= \inf_{\xi:\,
  \PPC \xi = \eta} \calR^*(\xi)$. 
\end{proposition}
\begin{proof} By Fenchel's equivalence we have $\xi\in \pl\Psi(y)\
  \Longleftrightarrow \ y \in \pl\Psi^*(\xi)$. Hence, the NESS equation
  \eqref{eq:NESSeqn} can be rewritten as 
\[
\rmD \calR \big(\ol u_\eta, \PPO  \ol y_\eta\big)= - \rmD\calE(\ol u_\eta), \quad \PPC
 \rmD\calE(\ol u_\eta ) =-\eta, \quad \ol y_\eta \in Y.
\]
Hence, we have the relation
$\eta = \PPC \rmD \calR \big( \ol u_\eta, \PPO \ol y_\eta\big) = \rmD_y
\sfR_Y(\ol u_\eta,y)$. Only if $\calR$ and hence $\sfR_Y$ are independent
of $u$, we obtain an explicit relation between $\eta$ and $\ol y_\eta$. 
Applying Fenchel's equivalence once again, we arrive at the assertion
$\ol y_\eta= \rmD_\eta \sfR_Y^*(\eta)$.

The second characterization of $\sfR_Y^*$ follows by an
application of Lemma \ref{le:MaaMie20}. 
\end{proof}

\subsection{Classical saddle points}
\label{su:ClassSP}

For a gradient system $(X,\calE,\calR)$ we consider the B-function 
\begin{equation}
  \label{eq:def.calL.ER}
  \frakB_{\calE,\calR}(u,\xi)= \calR^*(u,\xi) - \calR^*(u,-\rmD\calE(u))  ,
\end{equation}
which is defined on $X\ti X^*$. It will be the source of a series of results
concerning NESS. We will simply write $\frakB$ in place of $\frakB_{\calE,\calR}$
if the relevant GS $(X,\calE,\calR)$ is clear.

\begin{remark}[Slope dissipation term]
  In the definition of $\frakB_{\calE,\calR}$, we use the formula
  $\calR^*(u,-\rmD\calE(u))$ to denote the so-called $\calR$-\emph{slope}, 
  which should properly be defined by its weak lower semicontinuous hull,
  namely
\begin{equation}
  \label{eq:Slope.Hull}
  \calS_\calR(u):= \inf\Bigset{\liminf_{n\to \infty}\calR^*(u_n,{-}\rmD\calE(u_n)) }{
  u_n \weak u,\ u_n \in \dom(\rmD\calE) }. 
\end{equation}
For example the linear diffusion equation $\dot u=\Delta u$ with no-flux
boundary conditions is the gradient-flow equation associated with the Otto
gradient system $(\rmP(\Omega), \calE_{\rmB\rmz},\calR_\text{Otto}^*)$ with
$\calE_{\rmB\rmz}(u)=\int_\Omega \LB(u)\dd x$, and
$ \calR_\text{Otto}^*(u,\xi) = \int_\Omega \frac12 |\nabla \xi|^2 u\dd x$. We
obtain the Fisher information
$\calS_\calR(u)= \int_\Omega 2 |\nabla\sqrt{u}|^2 \dd x$, which is well defined even
when $u=0$ in a set of positive measure, whereas $u\in \dom(\rmD\calE)$ needs
$u>0$ a.e.

Subsequently, we will still write $\calR^*(u,-\rmD\calE(u))$ to emphasize the
structure of the problem; but whenever analysis is done, one has to replace
this term by $\calS_\calR$.
\end{remark}

Obviously, for all $ u \in X $ the functions $\frakB(u,\cdot):X^*\to \R$ are
convex, and in some cases we have concavity of $\frakB(\cdot,\xi)$ for all
$\xi\in X^*$. In the case of quadratic energy $\calE_\rmQ(u)=\frac12\langle \bbA
u,u\rangle - \langle \ell, u\rangle$ and a quadratic dual dissipation potential
$\calR^*_\rmQ(u,\xi)=\frac12\langle \xi, \bbK \xi\rangle$ we obtain the
simple quadratic B-function
\begin{equation}
\label{eq:calLquadr} 
\frakB_{\calE_\rmQ,\calR_\rmQ}(u,\xi)= 
\scrB_\mafo{quad}(u,\xi) = \frac12\langle \xi,\bbK\xi\rangle - \frac12
\big\langle \bbA u{-}\ell,\bbK (\bbA u{-}\ell) \big\rangle.  
\end{equation}
which has the above-mentioned concave-convex property on $X\ti X^*$. However,
our theory does not need concavity in $u$.

\begin{definition}[Global saddle points]
\label{de:GlobalSP} Given two Banach spaces $X$ and $Y$ and a functional
$\scrB:X\ti Y\to \R$, we call a point $(\ol x, \ol y)\in X\ti Y$ a \emph{(global)
  saddle point} for $\scrB$ if 
\[
{ \forall\, x\in X, \ y\in Y: \quad \scrB(x,\ol y) \leq \scrB(\ol x, \ol y) \leq
\scrB(\ol x,y). }
\]
\end{definition}

Thus, we are in the situation of classical saddle-point theory, see
\cite{EkeTem74ACPV} and Appendix \ref{se:Saddles}, which collects the most
important facts. In particular, we will use that the infimum over
$\xi\in X$ and the supremum over $u\in X$ can be interchanged if a saddle point
exists, see Lemma~\ref{le:SimpleFactsSP}:
\begin{subequations}
\label{eq:EstimSP}
\begin{align}
\text{(a)}\qquad &  \label{eq:EstimSP.a}
 \SI_{\scrB} := \sup_{u\in X} \inf_{\xi\in X^*} \scrB(u,\xi) \ \leq \ 
   \inf_{\xi\in X^*} \sup_{u\in X} \scrB(u,\xi) := \IS_{\scrB} , 
\\
\text{(b)}\qquad & \label{eq:EstimSP.b}
\text{saddle point $(\ol u, \ol \xi)$ exists} \ \ \Longrightarrow \ \ 
\SI_{\scrB}  =\IS_{\scrB}   = \scrB(\ol u, \ol \xi). 
\end{align}
\end{subequations}

For $\scrB_\mafo{quad}$ in \eqref{eq:calLquadr} with invertible $\bbA$ we see
that $(\ol u, \ol \xi)$ is a saddle point if and only if $\ol\xi=0$ (use
$\bbK>0$) and $ \rmD \calE (\ol u) = \bbA \ol u {-} \ell =0$, viz.\
$\ol u=\bbA^{-1}\ell$. We then have $\scrB_\mafo{quad}(\ol u, 0)=0$. If $\bbA$ is not
invertible, we have multiple saddle points, namely all $\ol u$ minimizing
$u\mapsto \frac12 \big\langle \bbA u{-}\ell,\bbK (\bbA u{-}\ell)
\big\rangle$. Then, one has 
$\scrB_\mafo{quad}(\ol u,0)= - \min\bigset{ \frac12 \langle \bbA u{-}\ell,\bbK
  (\bbA u{-}\ell) \rangle}{u \in X}$.

As a second example we consider 
\begin{equation}
  \label{eq:ExaPlast}
  X=\R^{i_*}, \quad \calE(u)=\frac12\langle \bbA u,u\rangle, \quad\text{and } \
\calR(v)=\sum_{i=1}^{i_*}\big( \sigma|v_i|+
\frac\nu2|v_i|^2\big)=\sigma|v|_1+\frac\nu2|v|_2^2. 
\end{equation}
Now we have $\calR^*(\xi)=\sum_{i=1}^{i_*} \frac1{2\nu} \big(
\max\{|\xi_i|{-}\sigma,0\}\big)^2$ such that $\calR^*(u,\xi)=0$ for
$|\xi|_\infty \leq \sigma$. Hence, we have many saddle points $(\ol u,
\ol\xi)$, namely all pairs with $|\ol\xi|_\infty\leq \sigma$ and $|\bbA\ol
u|_\infty \leq \sigma$. Again all saddle points satisfy $\frakB (\ol u,\ol
\xi)=0$.     

In the following, we give a general characterization of saddle points and 
complement the result with a discussion of critical points $(\wt u,\wt \xi)$ of
$\scrB$ satisfying $\wt\xi=-\rmD \calE(\wt u)$. 

\begin{theorem}[Unconstrained saddle points] 
\label{th:SaddleNESS} 
Consider a GS $(X,\calE,\calR)$ and set $\scrB=\frakB_{\calE,\calR}$ as in
\eqref{eq:def.calL.ER}. 

(a) A pair $(\ol u,\ol\xi)$ is a (global) saddle point of $\scrB$ if and only
if 
\[
\calR^*(\ol u,\ol\xi) =0\quad \text{and} \quad \calR^*(\ol u,-\rmD\calE(\ol
u))= \min_{u\in X} \calR^*(u,-\rmD\calE(u)).
\]

(b) If there exists $u_*\in X$ with $\rmD \calE(u_*)=0$, then
all saddle points satisfy the relation 
$\calR^*\big(\ol u,-\rmD\calE(\ol u)\big)=0$, and hence
$\scrB(\ol u,\ol\xi)=0$.

(c) If in addition to the condition in (b), the dual dissipation potentials
$\calR^*(u,\cdot) :X^*\to \R$ are strictly convex, then all saddle points  $(
\ol u, \ol\xi )$ satisfy $\ol\xi =0$ and $\rmD\calE(\ol u)=0$. 
\end{theorem}
\begin{proof}
\STEP{Part (a).}  Minimizing $\scrB$ with respect to $\xi\in X^*$ and using 
use $0 =\calR^*(u,0)\leq \calR^*(u,\xi)$ yields
\[
\SI_{\scrB}  = \sup_{u\in X} \big({-}\calR^*(u,-\rmD\calE(u))\big)=:\ol S \ \leq 0. 
\]
Moreover, choosing $\xi=0$ we obtain an upper bound for $\IS_{\scrB}  $, namely 
$\IS_{\scrB}  \leq \ol S$. Thus, with \eqref{eq:EstimSP.a} we conclude
$\SI_{\scrB}  =\IS_{\scrB}   =\ol S$. 

Hence, we conclude that a saddle point $(\ol u, \ol\xi)$ must satisfy
$\calR^*(\ol u, -\rmD\calE(\ol u))=-\ol S$ and $\calR^*(\ol u, \ol\xi )=0$,
which is the desired result (a).

\STEP{Part (b).} We obtain $\ol S=0$ and the result follows.

\STEP{Part (c).} This is an immediate consequence  of the implication 
$\calR^*(u,\xi)=0 \ \Rightarrow\ \xi=0$ and of Part (b).
\end{proof} 

The following result will not be used in the sequel, but it gives a first
insight why the saddle-point theory for $\frakB_{\calE,\calR}$ is useful. The
point is that there is a certain redundancy in the Euler-Lagrange equation for
critical points $(\wt u, \wt\xi)$ of $\frakB_{\calE,\calR}$, when the critical
point satisfies $\wt\xi = -\rmD\calE(\wt u)$. 

\begin{lemma}[Euler-Lagrange equations if $\wt\xi = -\rmD\calE(\wt u)$]
\label{le:EulLag}
The pair $(\wt u,\wt\xi)= (\wt u, -\rmD\calE(\wt u)) $ $ \in X\ti X^*$ is
a critical point of $\scrB= \frakB_{\calE,\calR} $ if and only  if 
$\rmD_\xi\calR^*( \wt u,-\rmD\calE(\wt u))=0 \in X$, i.e.\ $\wt u$ is a steady state
for the gradient system $(X,\calE,\calR)$. 
\end{lemma}
\begin{proof}
We have $\rmD_\xi \scrB(u,\xi)[\wt\xi]=\rmD_\xi\calR^*(u,\xi)[\wt\xi]$ and
\[
\rmD_u\scrB(u,\xi)[\wt u]=\rmD_u\calR^*(u,\xi)[\wt u]
-\rmD_u\calR^*(u,-\rmD\calE(u))[\wt u] +
\rmD_\xi\calR^*(u,-\rmD\calE(u))\big[\rmD^2\calE(u)[\wt u]\big]. 
\]
Inserting $(u,\xi)= (\wt u,-\rmD\calE(\wt u)) $ we see a
cancellation and the two equations for a critical point reduce to 
\[
0=\rmD_\xi \scrB(u_*,\xi_*)[\wt\xi]=\rmD_\xi\calR^*(u_*,\xi_*)[\wt\xi], \quad 
 0=\rmD_u\scrB(u_*,\xi_*)[\wt u] = \rmD_\xi\calR^*(u_*,\xi_*)
\big[\rmD^2\calE(u_*)[\wt u]\big]  .
\]
Thus, we see that it is necessary and sufficient to satisfy
$\pl_\xi\calR^*(u_*,\xi_*)=0$. 
\end{proof}

\begin{remark}[Perturbed gradient systems]\slshape
If a the gradient-flow equation is perturbed by a
general vector field $V$ in the form 
\begin{equation}
  \label{eq:PertGradFlow}
  \dot u = V(u) + \pl_\xi\calR^*(u,{-}\rmD \calE(u)),
\end{equation}
then steady states can still be obtained as stationary points of a B-function
$\wt\scrB$, namely 
\begin{equation*}
%  \label{eq:def.wt.mfL}
\wt{\scrB}(u,\xi)= \calR^*(u,\xi) -\langle \xi,V(u)\rangle
-\calR^*(u,{-}\rmD\calE(u)) - \langle \rmD\calE(u),V(u)\rangle.
\end{equation*}
Assume that  $u_*$ is a steady state for \eqref{eq:PertGradFlow}, namely 
\begin{equation}
  \label{eq:SteadyState}
  0=V(u_*) + \pl_\xi\calR^*(u_*,{-}\rmD\calE(u_*)),
\end{equation} 
then $(u,\xi)=(u_*,{-}\rmD\calE(u_*))$ is a stationary point for
$\wt{\scrB}$ and obviously the critical value is $0$, 
i.e.\ $\wt{\scrB}(u_*,{-}\rmD\calE(u_*))=0$.

To see the stationarity we observe 
$\rmD_\xi \wt\scrB(u,\xi)= \rmD_\xi \calR^*(u,\xi) - V(q)$,
and \eqref{eq:SteadyState} yields $
\rmD_\xi \wt\scrB(u_*,{-}\rmD\calE(u_*))=0$ as desired. 
For the derivative with respect to $u$ we have 
\begin{align*}
& \rmD_u\wt\scrB(u,\xi)[w]= \rmD_u\calR^*(u,\xi)[w] - \langle
\xi,\rmD V(u)[w]\rangle   -\rmD_u \calR^*(u,{-}\rmD \calE(u))[w]\\
&\qquad 
-\rmD_\xi \calR^*(u,{-}\rmD\calE(u))[{-}\rmD^2\calE(u)[w,\cdot]\rangle 
%\\ &\quad 
- \rmD^2\calE(u)[w,V(u)] - \rmD \calE(u)\big[ \rmD V(u)[w]\big].  
\end{align*}
Inserting $\xi = {-}\rmD\calE(u)$ the first term cancels
the third, and the second term cancels the last. Moreover, the forth
and the fifth terms cancel if we additionally use
\eqref{eq:SteadyState}. Hence,  $\rmD_u \wt\scrB(u_*,{-}\rmD\calE(u_*))=0$,
and $(u_*,{-}\rmD\calE(u_*))$ is indeed a stationary point for $\wt\scrB$. 
\end{remark}

\subsection{Constrained saddle points} 
\label{su:ConstrainedSP}

We now study the constrained case generated by a port GS
$(X,\calE,\calR,\PPO)$, where the linear port mapping $\PPO: Y \to X$
is used to drive the GS $(X,\calE,\calR)$. We start by introducing a
constrained saddle-point problem and then relate the existence of constrained
saddle points to the existence of NESS.

\begin{problem}[Constrained saddle-point problem (CSPP)]
\label{problem:CSPP}
Given the port GS $(X,\calE,\calR,\PPO)$ with B-function
$\frakB_{\calE,\calR}$ and port mapping $\PPO: Y \to X$, the
\emph{constrained saddle-point problem} for $\eta\in Y^*$ consists in finding a
saddle point $(\ol u_\eta,\ol\xi_\eta) \in X\ti X^*$ for
\begin{equation}
 \label{eq:CSPP}
\begin{aligned} 
\forall \; u \in X \text{ with } \PPC \rmD\calE(u)=-\eta& \ \ 
\forall \; \xi \in X^*\text{ with } \PPC \xi=\eta: 
\\  
&
\frakB_{\calE,\calR}(u,\ol\xi_\eta ) \leq \frakB_{\calE,\calR}(\ol u_\eta,\ol\xi_\eta ) \leq 
\frakB_{\calE,\calR}(\ol u_\eta,\xi). 
\end{aligned}
\end{equation}
The saddle point $ (\ol u_\eta,\ol\xi_\eta) \in X\ti X^*$ is called a
\emph{null-saddle} if $\frakB_{\calE,\calR}(\ol
u_\eta,\ol\xi_\eta )=0$.
\end{problem}

In light of our theory, it will be important to know whether a constrained
saddle point $(\ol u_\eta,\ol\xi_\eta)$ gives rise to a NESS. For
this we need the following preliminary result. 

\begin{lemma}[Null-saddles]
\label{le:NullSP.NESS}
If a constrained saddle point $ (\ol u_\eta,\ol\xi_\eta) \in X\ti X^*$ 
satisfies $\ol\xi_\eta= -\rmD\calE(\ol u_\eta)$, then it is a
null-saddle. Vice versa, if $\calR^*(\ol u_\eta ,\cdot):X^*\to [0,\infty]$ is
strictly convex and $ (\ol u_\eta,\ol\xi_\eta) \in X\ti X^*$ is a null-saddle,
then it satisfies $ \ol\xi_\eta = - \rmD\calE(\ol u_\eta)$.
\end{lemma}
\begin{proof}
  The first statement follows directly from
  $\frakB_{\calE,\calR}(u,-\rmD\calE(u))=0$ for all $u\in X$.

For the opposite implication we start from a null-saddle $(\ol
u,\ol\xi)$. From $0= \frakB_{\calE,\calR} ( \ol u,\ol\xi) \leq
\frakB_{\calE,\calR} (\ol u, \xi)$ for all $
\xi$ with $\PPC \xi=\eta$ we see that $\xi = \ol \xi$ and $\xi = -
\rmD\calE(\ol u)$ are global minimizers. By strict convexity the minimizer is 
unique, which proves the assertion.
\end{proof}

We recall the example in \eqref{eq:ExaPlast} where $\calR^*$ is not strictly
convex, because of $\calR^*(\xi)=0$ for $|\xi|_\infty\leq \sigma$. The saddle
points $(\ol u,\ol\xi)$ are characterized by $|\bbA \ol u|_\infty\leq \sigma$
and $|\ol\xi|_\infty\leq \sigma$ and all of them are null-saddles. However,
only some satisfy $\ol\xi = -\bbA\ol u$. This shows that the result does not
hold without a further condition like our strict convexity.

The next result shows that constrained saddle points of the form
$(\ol u_\eta,-\rmD\calE(\ol u_\eta))$ satisfies a simplified Euler-Lagrange
equation, and hence are NESS. Of course, without strong further global
assumptions one cannot expect that all NESS can be obtained as constrained
saddle points. 

\begin{proposition}[Euler-Lagrange equations for NESS] 
If the constrained saddle point in \eqref{problem:CSPP} has the form
$(\ol u,\ol \xi)=(\ol u_\eta,-\rmD\calE(\ol u_\eta))$, then $\ol u_\eta$
is a NESS for the port GS $(X,\calE,\calR,\PPO)$ with the constraint $\PPC
\rmD\calE(u)=-\eta$. 
\end{proposition}
\begin{proof}
In \eqref{eq:CSPP} we may consider variations $\wh\xi$ and $\wh u$ with $\PPC \wh
\xi=0$ and $\PPC\rmD^2\calE(u)[\wh u]=0$. Thus, we obtain
\begin{align*}
0&= \rmD_\xi \frakB(u,\xi)[\wh \xi]= \rmD_\xi\calR^*(u,\xi)[\wh\xi] = \langle
\wh\xi, \rmD_\xi\calR^*(u,\xi)\rangle_X,
\\
0&= \rmD_u\frakB(u,\xi)[\wh u] = \rmD_u \calR^*(u,\xi)[\wh u] - \rmD_u
\calR^*(u,-\rmD\calE(u))[\wh u] +
\rmD_\xi\calR^*(u,{-}\rmD\calE(u))\big[\rmD^2\calE(u)[\wh u]\big]. 
\end{align*}
Inserting $\ol\xi= {-}\rmD\calE(\ol u)$ we obtain a cancellation in the second
line leading to 
\[
0= \langle
\wh\xi, \rmD_\xi\calR^*(\ol u,{-}\rmD\calE(\ol u))\rangle_X \quad \text{and} \quad 
0= \langle \rmD^2\calE(\ol u)[\wh u],\rmD_\xi\calR^*(\ol u,\ol\xi)\rangle_X.
\]
However, by the choice of admissible variations, we see that the second
relation follows from the first. Hence we have
$\rmD_\xi\calR^*(\ol u,{-}\rmD\calE(\ol u)) \in
\big(\mafo{ker}(\PPC)\big)^\perp $.

To conclude, we simply use Fredholm's
alternative (theorem): 
\[
\big(\mafo{ker}(\PPC)\big)^\perp := \bigset{x\in X}{ \PPC\xi=0 \ 
   \Rightarrow\ \langle   \xi,x\rangle_X=0 } \;
 = \; \mafo{ran}(\PPO ):=\bigset{\PPO y}{ y \in Y}.
\]
With this we have
$\rmD_\xi\calR^*(\ol u,{-}\rmD\calE(\ol u)) \in \bigset{\PPO y}{ y \in
  Y}$, which gives $\ol y \in Y$ such that \eqref{eq:NESSeqn} holds.
\end{proof}

We now provide a general existence result for constrained saddle points and for
NESS. For this we use the following additional assumptions on
$\frakB=\frakB_{\calE,\calR}$ and $\PPO: Y\to X$: 
\begin{subequations}
\label{eq:CondsERP}
\begin{align}
&\label{eq:Conds.calL} 
 \left. \begin{aligned} 
\forall\, u \in X,\ \xi \in X^*:& \quad 
  \frakB(u,\cdot): X^*\to \R \ \text{ and } \ {-}\frakB(\cdot,\xi):X \to
\R \ \text{ are} 
\\
&\text{lower semi-continuous, strictly convex, and coercive};
\end{aligned}\right\}
\\
\label{eq:Conds.B}
&\forall \, \eta \in Y^*: \quad \bigset{u\in X}{ \PPC\rmD\calE(u)=-\eta}  
\text{ is nonempty, closed, and convex}. 
\end{align}
\end{subequations}

\begin{theorem}[Existence of constrained saddle points]
  \label{th:ExiCSP} Assume that $(X,\calE,\calR,\PPO)$ satisfies
  \eqref{eq:CondsERP}.  Then, for each $\eta \in Y^*$ there exists a unique
  constrained saddle point $(\ol u_\eta, \ol\xi_\eta)$ for $\frakB$ (in the
  sense of \eqref{eq:CSPP}).

If additionally the mapping $X\ni u \mapsto \rmD\calE(u) \in X^*$ is
surjective, then these saddles points are NESS satisfying
$\ol\xi_\eta=-\rmD\calE( \ol u_\eta)$
and \eqref{eq:NESSeqn}.   
\end{theorem}
\begin{proof} The existence follows by applying Proposition
\ref{pr:ExiSaddlePoint} with $\bfU=\bigset{\xi\in X^*}{ \PPC \xi = \eta}$ and
$\wt\bfV= \bigset{u\in X}{  \PPC \rmD\calE(u)=-\eta }$, where we extend $\frakB$ by
$-\infty$ outside of $\wt\bfV$ if it is not a linear space.  Thus, we find a
unique constrained saddle point $(\ol u_\eta, \ol\xi_\eta)$ with $\PPC
\rmD\calE( \ol u_\eta)= -\eta$ and $\PPC \ol\xi_\eta = \eta$.

Using Lemma \ref{le:NullSP.NESS} it is sufficient to show that $(\ol
u_\eta, \ol\xi_\eta)$ is a null-saddle. Because we already have a saddle point,
it is sufficient to show $\SI_{\frakB}   \leq 0 \leq \IS_{\frakB}  $. 

For the lower estimate we simply use $ \inf_{\xi\in \bfU} \frakB(u,\xi) \leq 
\frakB(u,-\rmD\calE(u))=0$. Taking the supremum over $u\in \wt\bfV$ we find
$SI_{\frakB}  \leq 0$. 

For the upper estimate we start from a general $\xi \in \bfU$ such that the
surjectivity of $\rmD\calE$ provides a $u_\xi\in \bfV$ with
$\xi=-\rmD\calE(u_\xi)$. With this we have
$\sup_{u\in \bfV} \frakB(u,\xi)\geq \frakB(u_\xi,\xi)=0$. Now taking the infimum
over $\xi \in \bfU$ yields $\IS_{\frakB} \geq 0$ as desired. 
\end{proof}

\subsection{NESS as minimizers}
\label{su:NESS.minimizer}

The main observation of the last section is that the equation
\eqref{eq:NESSeqn} does not have a simple variational structure. Its
characterization via the above saddle-point theory provides some kind of
variational structure, but needs a doubling of variables. Moreover, in
nonlinear problems (non-quadratic $\frakB$) the saddle-point theory for solving
infinite-dimensional problem like PDEs is technically rather demanding.

The naive way of treating the CSPP
\eqref{eq:CSPP} would be to minimize first with
respect to $\xi$ providing $\xi= \Xi_B(z,u)$ and such that it remains to study
the minimization problem
\[
u \ \mapsto \ \calR(u,{-}\rmD\calE(u)) - \calR^*(u,\Xi_B(z,u)) \quad
\text{subject to } B\rmD\calE(u)=z.
\]
This approach is doable but has the disadvantage that it is difficult to keep
enough control on the mapping $u \mapsto \Xi_B(z,u)$ to tackle the final
minimization problem.

The following result shows that the saddle point can be turned into a
minimization problem by applying a suitable Legendre transformation with
respect to the constrained variable $\xi$, but keeping a dual parameter
$\Lambda \in Z^*$. Thus, the minimization formulation stays explicit in terms
of the constituents of the GS $(X,\calE,\calR)$. Moreover, it is more
directly related to the Euler-Lagrange equations \eqref{eq:NESSeqn} and
the original thermodynamical functions.

\begin{proposition}[NESS as minimizers]
\label{pr:NESS.Min.Lag}
For all $\eta\in Y^*$ any global minimizer $(\ol u, \ol y) \in X\ti Y$ of the
constrained minimization problem
\begin{equation}
  \label{eq:Min.NESS}
\begin{aligned}
\text{minimize } \ &\calR(u,\PPO y) + \calR^*(u,-\rmD\calE(u)) +  \langle
\eta,y\rangle_Y\\
\text{over }\ & (u,y) \in X\ti Y \quad 
\text{subject to } \ \PPC \rmD\calE(u)= -\eta
\end{aligned}   
\end{equation}
gives rise to a constrained saddle points $(\ol u,\ol\xi) \in X\ti X^*$ for
\eqref{eq:CSPP} where we can choose any $\ol\xi \in 
 \mafo{Argmin}\bigset{\calR^*(\ol u, \xi)}{  \PPC \xi=\eta}$. Vice versa,
if $(\ol u,\ol\xi)$ is a constrained saddle point for \eqref{eq:CSPP}, then
$(\ol u,\ol y)$ with $\ol y \in  \mafo{Argmax}\bigset{\langle \eta,y\rangle -
  \calR(\ol u, \PPO y)}{ y\in Y}$ is a global minimizer for \eqref{eq:Min.NESS}.

Moreover, if $(\ol u,\ol y)$ is a null-minimizer, then $(\ol y,\ol\xi)$ is a
null-saddle, and under the additional assumption of strict convexity of
$\calR^*(\ol u,\cdot)$ it defines a NESS solving \eqref{eq:NESSeqn}. 
\end{proposition}
\begin{proof} 
  We define the auxiliary dissipation potentials
  $\Psi_u:Y\to \R_\infty;\, y \mapsto \calR(u,\PPO y)$ and can now apply
  Lemma \ref{le:MaaMie20} below. This gives
\begin{align}
  \label{eq:MaaMieRelation}
\inf_{\rule{0pt}{0.6em}\xi:\PPC \xi=\eta} \calR^*(u,\xi)& 
  = \Psi^*_u (\eta)  = \sup_{y \in Y} 
 \big( \langle \eta, y\rangle_Y {-} \calR(u,\PPO  y)   \big). 
\end{align}

With this we obtain the following chain of identities:
\begin{align*}
& \mbox{}\hspace{-1em}\sup_{u\in X \atop \PPC\rmD\calE(u)=-\eta} \ 
\inf_{\xi \in X^* \atop   \PPC \xi=\eta}
\frakB_{\calE,\calR}(u,\xi)
= \sup_{u\in X \atop \PPC\rmD\calE(u)=-\eta}
 \Big( \big[\inf_{\xi \in X^* \atop   \PPC \xi=\eta} 
\calR^*(u,\xi) \big] - \calR^*(u,{-}\rmD \calE(u)) \Big) 
\\
&\overset{\text{\eqref{eq:MaaMieRelation}}}= 
\sup_{u\in X \atop \PPC\rmD\calE(u)=-\eta}  \Big(\  \Big[ 
\sup_{y\in Y} 
 \big( \langle \eta, y \rangle_Y -\calR(u,\PPO  y) \big)
 \Big] - \calR^*(u,{-}\rmD \calE(u)) \Big) 
\\
& = - \inf_{u:\PPC \rmD\calE(u)=-\eta \atop y \in Y} \Big( \calR(u,\PPO  y)
+ \langle \eta, y\rangle_Y  + \calR^*(u,{-}\rmD \calE(u)) \Big).   
\end{align*}
This shows that the minimization problem \eqref{eq:Min.NESS} is equivalent to
the CSPP \eqref{eq:CSPP} if we choose $\xi=\ol\xi \in X^*$ in
\eqref{eq:MaaMieRelation} optimally, i.e.\
$\ol\xi =\mafo{Argmin}\bigset{\calR^*(\ol u, \xi)}{ \PPC \xi=\eta}$.

Moreover, the values are the same up to a minus sign. Hence, null-minimizers
$(\ol u,\ol y) \in X\ti Y$ correspond to null-saddles  $(\ol u,\ol\xi) \in X
\ti X^*$, and the remaining statement follows from Lemma
\ref{le:NullSP.NESS}. 
\end{proof}

In the above proof the relation in \eqref{eq:MaaMieRelation} relies on the
following general result.

\begin{lemma}
\label{le:MaaMie20} 
For a lower semi-continuous and convex $\Psi:X\to \R_\infty$ and linear bounded
operator $B:X^*\to Z$  we have 
\[
\inf_{\rule{0pt}{0.4em}\xi \in X^*:\; B\xi=z} \Psi^*(\xi)= \sup_{\Lambda \in Z^*}
\Big( \langle \Lambda, z \rangle_Z - \Psi(B^*\Lambda) \Big) . 
\]  
\end{lemma}
\begin{proof} Consider a dissipation potential $\Psi:X\to [0,\infty]$ and a
  bounded linear mapping $A: Y\to X$ and define the dissipation potential
  $\wt\Psi:Y\to [0,\infty];\ y \mapsto \Psi(A y)$. In
  \cite[Prop.\,6.1]{MaaMie20MCRS} the identity $ (\wt\Psi){}^*(\eta)
  = \inf\bigset{\Psi^*(\xi)}{ A^*\xi = \eta}$ is established.
Applying this with $Y=Z^*$ and $A=B^*:Z^*\to X$ the
assertion follows.  
\end{proof}

A much simpler case occurs if the dissipation potential $\calR$ is independent
of the state $u\in X$. Then, the minimization in \eqref{eq:Min.NESS} with
respect to $u\in X$ subject to $\PPC\rmD\calE(u)=-\eta$ and with respect to
$y \in Y$ decouple completely. In particular, if $\calE$ is uniformly convex,
all
\[
 \ol u \in \mafo{Arg\!\;min}\bigset{\calR^*({-}\rmD\calE(u)) }{
   \PPC\rmD\calE(u)=-\eta}
\]
are NESS, see \cite{Miel23?PGSN} for more details. This relates to Prigogine's
principle that states that NESS are minimizers of the dissipation, i.e.\ of
$u\mapsto \Phi_*({-}\rmD\calE(u))$, where $\Phi_*\geq \calR^*$ is defined after
\eqref{eq:NESSeqn}. For quadratic dissipation potentials (linear kinetic
relations) we have $\Phi_*(\xi)=2\calR^*(\xi)$, such that the principle becomes
exact. For general nonlinear kinetic relations the result is an approximation
only, which works well close to equilibrium, see the discussions in 
\cite[Cha.\,V]{DegMaz84NET} and \cite{Miel23?PGSN}.

\subsection{Constrained B-functions, BER structure, and NESS} 
\label{su:ReducLagrDual}
 
When doing reduction or $\Gamma$-limits of B-functions, we may end up with a
general function $\scrK:Y\ti Y^* \to \R$ and may then ask the question whether
this function can be written as a B-function $\frakB_{\sfE,\sfR}$. 
 
\begin{definition}[BER structure]
\label{de:DualityStruct}
We say that a function $\scrK:Y\ti Y^* \to \R$ 
has the \emph{BER structure} $(\sfE,\sfR)$, if $(Y,\sfE,\sfR)$ is a
gradient system and  
\[
{\scrK=\frakB_{\sfE,\sfR}, \quad \text{namely } \forall\, (y,\eta)\in Y \ti
Y^*:\ \scrK(y,\eta)=\sfR^*(y,\eta) - \sfR^*(y,-\rmD\sfE(y)).}
\]
\end{definition} 

We observe that for a given $\scrK$ the dissipation functional $\sfR$ and its
dual $\sfR^*$ are uniquely determined by
$\sfR^*(y,\eta)=\scrK(y,\eta)-\scrK(y,0)$. Hence, we have the following
necessary and sufficient conditions of a BER structure. The third
condition \eqref{eq:Conds.DualStru.c} provides the important link to
null-saddles. 

\begin{proposition}[Conditions for  BER structure]
\label{pr:CondsDualStruc}
Given  an energy $\sfE:Y\to \R$, the function $\scrK:Y\ti Y^* \to \R$ has a BER
structure $(\sfE,\sfR)$ if and only if
\begin{subequations}
  \label{eq:Conds.DualStru}
\begin{align}
  \label{eq:Conds.DualStru.a}
& \forall\, (y,\eta)\in Y\ti Y^*:\quad \scrK(y,\eta)\geq \scrK(y,0),
\\
  \label{eq:Conds.DualStru.b}
& \forall\, y \in Y:\quad \scrK(y,\cdot):Y^*\to \R \ \text{ is convex}, 
\\
  \label{eq:Conds.DualStru.c}
&  \forall\, y \in Y:\quad \scrK(y,{-}\rmD \sfE(y))=0.
\end{align}
\end{subequations}
Then, $\sfR$ is given by $\sfR^*(y,\eta)=\scrK(y,\eta)-\scrK(y,0)$. 
\end{proposition}
\begin{proof} It is obvious that $\scrK$ satisfies \eqref{eq:Conds.DualStru} if
it has the BER structure $(\sfE,\sfR)$. 

To show the opposite, we observe that
$\sfR^*_\scrK: (y,\eta)\mapsto \scrK(y,\eta){-}\scrK(y,0)$ is a (dual)
dissipation potential because of \eqref{eq:Conds.DualStru.a} and
\eqref{eq:Conds.DualStru.b}. Inserting the formula for $\sfR^*_\scrK$ into the
condition $0=\frakB_{\sfE,\sfR}(y,\eta) -\sfR^*(y,\eta) +
\sfR^*(y,-\rmD\sfE(y))$ defining BER structures, we obtain
\begin{align*}
0&= \scrK(y,\eta)-\sfR^*_\scrK(y,\eta) + \sfR^*_\scrK(y,-\rmD\sfE(y))
\\
&= \scrK(y,\eta) -\big( \scrK(y,\eta){-}\scrK(y,0)\big) +
\big(\scrK(y,-\rmD\sfE(y))-\scrK(y,0) \big)   =  \scrK(y,-\rmD\sfE(y)).
\end{align*}
Hence, \eqref{eq:Conds.DualStru.c} guarantees that this $(\sfE,\sfR_\scrK)$ is
the desired BER structure.
\end{proof}

We return to our constrained saddle-point problems by generalizing it in a
crucial way. For this we use a second port function $P:X\to Y$ which allows
us to impose direct conditions $Pu=y$ on the state variable, whereas
$\PPC \rmD\calE(u)=\eta$ does this indirectly. Nevertheless, we always
assume there is an energy $\sfE:Y\to \R$, such that
\begin{equation}
  \label{eq:PPcircEy}
Pu = y \quad \Longrightarrow \quad \PPC \rmD\calE(u)=\rmD \sfE(y).
\end{equation}

An important point for understanding the reduced or effective kinetic
relation generated by the port GS $(X,\calE,\calR,\PPO)$ is to study the
reduced B-function $\scrB_\red:Y\ti Y^*\to \R$ defined via
\begin{equation}
  \label{eq:def.Lred}
  \scrB_\red(y,\eta):= \sup_{u\in X\atop Pu=y} \inf_{\xi\in X^*\atop \PPC
  \xi=\eta} \frakB_{\calE,\calR}(u,\xi) .
\end{equation}
In contrast to the previous analysis, we are now using two independent
constraints $y\in Y$ and $\eta\in Y^*$, whereas in Section
\ref{su:ConstrainedSP} we always assumed the compatibility
$\eta= - \rmD \sfE(y)$, cf.\ \eqref{eq:PPcircEy}. 
However, assuming there are null-saddles under these constraints means that
$\scrB_\red(y,- \rmD \sfE(y) )=0$ holds, i.e.\ the necessary 
\eqref{eq:Conds.DualStru.c} holds. The next result provides the
fundamental link between null-saddles and a BER structure for
$\scrB_\red$. 

\begin{theorem}[BER structure for $\scrB_\red$]
\label{th:Lred.DualityStr} 
Consider a gradient system $(X,\calE,\calR)$ with port mappings $P:X\to Y$ and
$\PPC:X^*\to Y^*$ and a compatible energy $\sfE$ as in \eqref{eq:PPcircEy}.
Assume that for all $y\in Y$ the CSPP \eqref{eq:CSPP} with
$\eta = -\rmD\sfE(y) = - \PPC \rmD\calE(u)$ has a null-saddle. Then, the
reduced B-function $\scrB_\red$ defined in \eqref{eq:def.Lred} has
the BER structure $(\sfR,\sfE)$ with
$\sfR^*(y,\eta)= \scrB_\red(y,\eta) -\scrB_\red(y,0)$, namely
\begin{equation}
  \label{eq:Lred.sfR*}
  \scrB_\red(y,\eta) { = \frakB_{\sfE,\sfR}(y,\eta) }
   = \sfR^*(y,\eta) -\sfR^*\big(y,-\rmD\sfE(y)\big).
\end{equation}
\end{theorem}
\begin{proof} The proof follows by checking the conditions
  \eqref{eq:Conds.DualStru} in Proposition \ref{pr:CondsDualStruc}. 

\STEP{Part (a):} Since $\calR^*$ is a dual dissipation potential we have 
\begin{align*}
\inf_{\xi:\:\PPC\xi=0} \frakB(u,\xi) = \frakB(u,0) \leq
\inf_{\xi:\:\PPC\xi=\eta} \frakB(u,\xi) .
\end{align*}
Taking the supremum over $u$ with $Pu=y$ gives $\scrB_\red(y,0)\leq
\scrB_\red(y,\eta)$ as desired.\medskip

\STEP{Part (b):} Defining $\calN(u,\eta)=\inf_{\xi:\,\PPC \xi=\eta}
\frakB(u,\xi)$ we can easily check that each $\calN(u,\cdot)$ is still
convex, because $\PPC $ is linear map. Indeed, for $\eta_0,\eta_1 \in Y$ and $\theta \in [0,T]$ set
$\eta_\theta = (1{-}\theta)\eta_0 + \theta \eta_1$. For $\eps>0$ pick
$\xi_0,\xi_1 \in X^*$ with $\frakB(u,\xi_j)\leq \calN(u,y_j)+\eps$. Then,
\begin{align*}
\calN(u,\eta_\theta) &\leq \inf_{\xi:\,\PPC \xi=\eta_\theta} \frakB(u,\xi)
\leq \frakB\big(u,(1{-}\theta)\xi_0 + \theta \xi_1\big)\\
&\overset{\frakB(u,\cdot)\text{ cvx}}\leq (1{-}\theta)
\frakB\big(u,\xi_0)  + \theta\frakB\big(u,\xi_1) \leq (1{-}\theta)
\calN(u,\eta_0 ) + \theta \calN(u,\eta_1) + \eps.
\end{align*}
Since $\eps>0$ was arbitrary, the convexity of $\calN(u,\cdot) $ is
established. 

Because $\scrB_\red(y,\cdot)$ is the supremum of the family
$\big(\calN(u,\cdot)\big)_{u:\,Pu=y}$ it is again convex.\medskip

\STEP{Part (c):} For every $y \in Y$ there exists a null-saddle $(\ol u,\ol
\xi)$ with $P \ol u=y$ and $\PPC \ol\xi= -\PPC \rmD\calE(\ol u)=-
\rmD\sfE(y)$. Hence, we have 
\[
\sup_{u:\,Pu=y} \frakB(u,\ol\xi) \leq \frakB(\ol u,\ol\xi)=0 \leq \inf_{\xi:\,
  \PPC \xi=-\rmD \sfE(y) } \frakB(\ol u,\xi) . 
\]
Comparing with the definition of $\scrB_\red$ we find
$\scrB_\red(y,{-}\rmD\sfE(y))=0$ as desired. 
\end{proof}

\section{EDP-convergence for slow-fast GSs via NESS}
\label{se:ReductGSNESS}

We consider a family of Gs $(X,\calE_\eps,\calR_\eps)$ where
$\eps>0$ is the small parameter modeling the ratio between fast and slow
relaxation times.  We consider two distinguished cases: in the first the state
space can be decomposed in the form $u= (U,w)\in X_\slow\ti X_\fast =X$ and in
the second we have
\[
X=\bigset{u=(U,w)\in X_\slow\ti X_\fast} { \bbQ(U,w):=P_\slow U-P_\fast w=0}
\]
where $P_\slow: X_\slow \to Y$ and $P_\fast : X_\fast \to Y$ are suitable port
mappings. Here we consider $U\in
X_\slow$ as the slow macroscopic part of the state variables, while $w\in
X_\fast$ is the fast microscopic part, that one wants to eliminate in the limit
$\eps\to 0$.

In both setting we assume that the scaling in $\eps$ is very particular, but
nevertheless we are able to treat a number of prototypical cases. In
particular, we assume $\calE_\eps(U,w)=E(U)+ \eps\, e(w)$.

\subsection{Case 1: product space $X= X_\slow \ti X_\fast$}
\label{su:Case1Product}

The precise assumptions on the scaling with $\eps>0$ are the following:
\begin{subequations}
  \label{eq:Scaling1.calReps}
\begin{align}
  \label{eq:Scaling1.calReps.a}
\calE_\eps(U,w) &= E(U) + \eps \,e(w) && \text{additive split of energy},
\\
  \label{eq:Scaling1.calReps.b}
 \calR^*_\eps(U,w;\Xi,\mu) &= \ol\calR^* (U,w; \Xi, \tfrac{\ds1}{\ds\eps}
\mu\big) &&\text{fast relaxation of $w$},
\end{align}
\end{subequations}
where $\ol\calR^*:X\ti X^*\to [0, \infty]$ is a general dual dissipation
potential independent of $\eps$. 

The associated gradient-flow equation takes a simple form, because
the appearance of $\eps$ is chosen in a particular way. 
\begin{subequations}
  \label{eq:Case1.eps.Eqn}
\begin{align} 
  \label{eq:Case1.eps.Eqn.a}
\dot U& = \rmD_\Xi\ol\calR^*\big( U,w;-\rmD E(U),-\rmD e(w) \big),
\\
  \label{eq:Case1.eps.Eqn.b}
\eps \dot w &= \rmD_\mu\ol\calR^*\big( U,w;-\rmD E(U),-\rmD e(w) \big).
\end{align}
\end{subequations}
Thus, on the formal level, we can drop the term $\eps \dot w$, because $w$
relaxes into a NESS on the time scale $\eps$ which is much faster than the
evolution of $U$ which happens on time scales of order $1$. The microscopic
variable $w$ moves into the NESS $w=\wt\sfw(U)$ satisfying
\begin{equation}
  \label{eq:GaCvg.1.NESS}
  0 = \rmD_{\mu} \ol\calR^*\big( U,w;-\rmD E(U), -\rmD e(w) \big) .
\end{equation}

Inserting the limiting relation $w=\wt\sfw(U)$ into the first equation of 
\eqref{eq:Case1.eps.Eqn.a} we obtain the reduced macroscopic problem 
\begin{equation}
  \label{eq:Case1.ReducEqn}
  \dot U = \rmD_\Xi \ol\calR^* \big(U,\wt\sfw(U); -\rmD E(U), -\rmD 
    e(\wt\sfw(U))\big).  
\end{equation}
The disadvantage of the above approach is that we lose control over
the gradient structures. As we have started with the GSs
$(X,\calE_\eps,\calR_\eps)$, it is natural to ask whether the effective
equation \eqref{eq:Case1.ReducEqn} has a natural gradient structure inherited
from $E$, $e$, and $\ol\calR$. 
\medskip

This question can be answered by the notion of EDP-convergence, which provides
a tool to stay on the level of gradient systems. We follow here the approach
developed in \cite{LMPR17MOGG} which forms the basis of the further developments
of EDP-convergence in \cite{DoFrMi19GSWE, MiMoPe21EFED}. The abbreviation
``EDP'' stand for the \emph{energy-dissipation principle} (cf.\
\cite[Thm.\,3.3.1]{Miel16EGCG}) that shows that under suitable technical
assumptions a curve $u_\eps=(U_\eps,w_\eps):[0,T]\to X$ is a solution of the
gradient-flow equation \eqref{eq:Case1.eps.Eqn} if and only if it satisfies the
energy-dissipation inequality
\[
\calE_\eps(u_\eps(T)) +\int_0^T\!\!\Big( \calR_\eps \big( u_\eps;\dot
u_\eps\big) + \calR^*_\eps \big( u_\eps; -\rmD\calE_\eps(u_\eps)\big) 
 \Big) \dd t \leq  \calE_\eps(u_\eps(0)). 
\] 
The idea in \cite{LMPR17MOGG, MaaMie20MCRS} is to replace the primal
dissipation $\calR_\eps(u,\dot u)$ by the lower bound $\langle \xi,\dot
u\rangle - \calR^*_\eps(u,\xi)$ for an arbitrary test function $\xi:[0,T]\to
X^*$. Then, the limit $\eps \to 0$ is performed and finally one maximizes with
respect to $\xi$ to recover the limiting energy-dissipation balance again. 

Thus, for a general smooth function $\xi :[0,T] \to X^*$ we have 
\[
\calE_\eps(u_\eps(T))
 + \int_0^T \!\!\Big( \big\langle \xi, \dot u_\eps\big\rangle
 - \calR^*_\eps \big( u_\eps; \xi \big)  
+ \calR^*_\eps \big( u_\eps; -\rmD\calE_\eps(u_\eps)\big) \Big) \dd t \leq 
\calE_\eps(u_\eps(0)). 
\] 
Using the explicit $\eps$-dependence of $\calE_\eps$ and $\calR^*_\eps$ imposed in
\eqref{eq:Scaling1.calReps} and choosing $\xi = (\Xi, \eps \zeta)$  we arrive
at 
\[
\calE_\eps(u_\eps(T)) 
+ \int_0^T\!\!\Big( \big\langle (\Xi,\eps\zeta), \dot u_\eps \big\rangle
 - \ol\calR^*\big( u_\eps; \Xi,\zeta \big)  + 
 \ol\calR^* \big( u_\eps; -\rmD E(U_\eps),-\rmD e(w_\eps)\big) \Big) \dd t \leq 
\calE_\eps(u_\eps(0)).  
\] 
Now passing to the limit $\eps\to 0$ the term
$\langle \eps \zeta,\dot w_\eps\rangle $ and the terms $\eps e(w_\eps(t))$
vanish. Assuming $(U_\eps,w_\eps)\to (U,w)$ we arrive at the inequality
\begin{align*}
E(U(T)) +\!\int_0^T \!\!\!\Big( \big\langle \Xi, \dot
U\big\rangle - \frakB_{\ol\calE,\ol\calR}(U,w; \Xi,\zeta) \Big) \dd t
 \leq  E(U(0)) \ \text{ for all }(\Xi,\zeta)\in \rmL^\infty([0,T];X^*), 
\end{align*}
where $\ol\calE(U,w)=E(U){+} e(w)$ and hence 
\[
\frakB_{\ol\calE,\ol\calR}(U,w; \Xi,\zeta) = \ol\calR^*\big(U,w; \Xi,\zeta \big) 
- \ol\calR^* \big(U,w; -\rmD E(U),-\rmD e(w)\big). 
\]

Since $w$ appears in the integral only via $w(t)$, but not with a derivative
$\dot w(t)$ we can eliminate $w(t)$ by taking the infimum pointwise in
$t\in [0,T]$. Similar, we can
eliminate $\zeta$ by a pointwise supremum. Hence, defining
$\scrB_\red : X_\slow\ti X^*_\slow\to \R$ via 
\begin{equation}
  \label{eq:def.calLslow}
  \scrB_\red(U,\Xi):=  \sup_{w\in X_\fast} \inf_{\zeta \in X_\fast^*}\; \!
\frakB_{\ol\calE,\ol\calR}(U,w; \Xi,\zeta). 
\end{equation}
we obtain the inequality 
\begin{align}
\label{eq:Case1.EDI.Lred}
E(U(T)) +\int_0^T \!\!\Big( \big\langle \Xi, \dot U\big\rangle 
 -  \scrB_\red(U,\Xi) \Big) \dd t  \leq  E(U(0)). 
\end{align}
Now it remains to show that $ \scrB_\red $ has a BER structure
$(E,\calR_\eff)$ in the sense of Definition \ref{de:DualityStruct}, i.e.\ it
has the form
\begin{align}
  \label{eq:Case1.Lred.Dual}
 \scrB_\red(U,\Xi)= \calR^*_\eff (U;\Xi) - \calR^*_\eff (U;-\rmD E(U)),\ \
 \text{ i.e. } \scrB_\red = \frakB_{E,\calR_\eff}
\end{align}
for a suitable effective dissipation potential $\calR_\eff$. 

If this is the case, we can insert this into \eqref{eq:Case1.EDI.Lred} and
reverse the Legendre transform with respect to $\Xi$ to obtain the
energy-dissipation inequality 
\begin{align}
\label{eq:Case1.EDI.Lred2}
E(U(T)) +\int_0^T \!\!\Big(  \calR_\eff (U; \dot U) + \calR^*_\eff \big(
U;-\rmD E(U) \big) \Big)  \dd t  \leq  E(U(0)). 
\end{align}
Applying the energy-dissipation principle once again, we see that $U$ is a
solution of the gradient-flow equation
\[
\dot U = \rmD\calR^*_\eff (U, {-}\rmD E(U))
\]
for the reduced gradient system $(X_\slow,E,\calR_\eff)$. Clearly, this
equation must equal \eqref{eq:Case1.ReducEqn}, but now we have a truly
thermodynamical structure.\medskip

To achieve this goal it remains to establish the BER structure
\eqref{eq:Case1.Lred.Dual}. The following result is the analogue of Theorem
\ref{th:Lred.DualityStr}.

\begin{theorem}[$\scrB_\red$ has BER structure]
\label{th:Lred.Case1.Dual} 
For a GS $(X_\slow\ti X_\fast,\ol\calE,\ol\calR)$ with
$\ol\calE = E {\oplus}\!\; e $ 
define $\scrB_\red: X_\slow \ti X_\slow^*\to \R$ as in
\eqref{eq:def.calLslow}. If for all $U\in X_\slow$ we have that
\[
\scrB_\red (U,{-}\rmD E(U)):= \sup_{w\in X_\fast}\   \inf_{\zeta\in X^*_\fast}
\frakB_{\ol\calE,\ol\calR}(U,w;-\rmD E(U),\zeta)  
\]
is a null-saddle (i.e.\ $\scrB_\red \big(U,{-}\rmD E(U)\big)=0 $), then
$\scrB_\red $ has the BER structure $(E,\calR_\eff)$ where $\calR_\eff$ is
given via $\calR^*_\eff(U,\Xi)= \scrB_\red(U,\Xi) -\scrB_\red(U,0)$.
\end{theorem}
\begin{proof} The result follows via Theorem
  \ref{th:Lred.DualityStr} if we use $Y=X_\slow$ and the port mappings 
\[
P(U,w)=U \in Y \quad \text{and} \quad  \PPO y = (y,0)\in X_\slow \ti X_\fast.
\]
Note that $\ol\calE =E{\otimes} e$ satisfies
$\rmD \ol\calE (U,w)=\big(\rmD E(U), \rmD e(w)\big)$, hence, $E:X_\slow \to \R$
is a compatible energy in the sense of \eqref{eq:PPcircEy}.
\end{proof}

\subsection{Case 2: factored product space $X=(X_\fast\ti
  X_\slow)\big/_{\!\!\mafo{ker} \,\bbQ}$}
\label{su:Case2Product}

In some cases it is not easy to decompose the state space $X$ into a product
$X_\slow\ti X_\fast$, but it is possible to decompose the state with some
overlay or joint traces on an interface, namely
 \[
X=\bigset{u=(U,w)\in X_\slow\ti X_\fast} { \bbQ(U,w):=P_\slow U-P_\fast w=0}
\]
where $P_\slow: X_\slow \to Y$, $P_\fast : X_\fast \to Y$,
$\PPC_\slow: X^*_\slow \to Y^*$, and $\PPC_\fast : X^*_\fast \to Y^*$ are
suitable port mappings. Below we will show that the chosen ansatz applies to
diffusion problems, where $P_\slow$ and $P_\fast$ are used to define traces
from two different sides of an interface, see \eqref{eq:DiffMembCompat} in
Section \ref{su:DiffMembrPDE}.
 
The precise assumptions are the following:
\begin{subequations}
  \label{eq:Scaling.calReps}
\begin{align}
  \label{eq:Scaling.calReps.a}
\calE_\eps(U,w) &= E(U) + \eps e(w) && \text{additive split of energy},
\\
  \label{eq:Scaling.calReps.b}
 \calR^*_\eps(U,w;\Xi,\xi) &= \wt\calR^* (U,w; \Xi, \tfrac{\ds1}{\ds\eps}
\xi\big) &&\text{fast relaxation of $w$},
\\
\nonumber
 \wt\calR^*(U,w;\Xi,\zeta)& = \calR^*_\slow(U;\Xi ) + \calR^*_\fast(w;\zeta) 
\\   \label{eq:Scaling.calReps.c}
&\quad + \bm\delta_{\{0\}}\big( \PPC_\fast\zeta{-} \PPC_\slow \Xi) &&
\text{interaction through $Y^*$}.
\end{align}
\end{subequations}
In principle, we could allow the more general case $ \calR^*_\fast(U,w;\zeta) $
in place of $\calR^*_\fast(w;\zeta) $, but refrain from doing so, because the
restricted version better highlights the fact that $U$ and $w$ mainly
interact through the ports via $Y$.

\TTODO{Do we need $\sfE:Y\to \R$ such that 
\[
P_\slow U=y \Rightarrow \PPC \rmD E(U)=\rmD E(y) \quad \text{and} \quad 
P_\fast w=y \Rightarrow \PPC \rmD e(w)=\rmD E(y) \ ?????
\]} 
Here $\bm\delta_{\{0\}}: Y^0 \to [0,\infty]$ is the convex function with
$\bm\delta_{\{0\}}(0)=0$ and $\infty$ otherwise. This function implements the
constraint $ \PPO_\fast^*\zeta = \PPO_\slow^* \Xi $ giving the interaction
condition $ \PPO_\fast^*\rmD e(w)= \PPO_\slow^* \rmD E(U)$. The subdifferential
of $\bm\delta_{\{0\}}$ at $\eta=0$ is given by $\pl\bm\delta_{\{0\}}(0)=Y$,
i.e.\ the hard constraint can transmit the fluxes
$(-\PPO _\slow y, \PPO _\fast y)$ for arbitrary $y \in Y$.

As before, we first observe that the gradient-flow equation takes a simple
form, because the appearance of $\eps$ is chosen in a particular way.
\begin{align*} 
&\binom{\dot U}{\eps \dot w} \in \pl\ol\calR^*\big( U,w;-\rmD E(U),-\rmD e(w) \big) 
\ \Longleftrightarrow \ 
\\
&\binom{\dot U}{\eps \dot w}= 
\binom{\rmD_\Xi \calR^*_\slow \big( U,-\rmD  E(U)\big)} 
      {\rmD_{\zeta} \calR^*_\fast\big( w,-\rmD e(w) \big) }
 + \binom{ \PPO_\slow y}{ {-}\PPO_\fast y}  \quad \text{with }
\left\{\ba{c}  P_\slow U=P_\fast w\\ \text{and } y\in Y. \ea \right.
\end{align*}

Thus, on the formal level, we can drop the term $\eps \dot w$, because $w$
relaxes into a NESS on the time scale $\eps$ which is much faster than the
evolution of $U$ which happens on time scales of order $1$. The microscopic
variable $w$ moves along the family of NESS $w=\wt\sfw(y)$ generated by the
flux $y\in Y$ via
\[
0 = \rmD_{\zeta} \calR^*_\fast\big( U,w;-\rmD e(w) \big) + \PPO _\fast y. 
\]

\TTODO{IMPROVE THIS:
We may assume that for each $\eta \in Y^*$, there exists a unique NESS
$w=\wt w(\eta)$ for the equation 
\[
0=  \pl_{\zeta} \calR^*_\fast \big( w,-\rmD e(w)\big) + \PPO  v, \quad
\PPC \rmD e(w)=\eta.
\]
and define the port mapping $\Phi:Y^*\to Y; \eta \mapsto v$. The constraint is
mapped to the Lagrange multiplier. 

The above theory shows that $\Phi$ has the following particular structure,
namely $v = \pl\calR^*_\red(y,\eta)$.}

As in the previous subsection, we can now involve the energy-dissipation
principle to show EDP-convergence, where now $\ol\calR^*$ is replaced by
$\wt\calR^*$ containing the constraint $\PPC_\slow\Xi=\PPC_\fast \zeta$. 
We again arrive at the reduced energy inequality \eqref{eq:Case1.EDI.Lred},
where now $\scrB_\red $ is replaced by $\scrB_\eff$ which takes a
special form due to the additive splitting of $\wt\calR^*$ in 
\eqref{eq:Scaling.calReps.c}: 
\begin{align}
  \nonumber
\scrB_\eff(U,\Xi) &=  \sup_{w\in X_\fast} \inf_{\zeta \in X_\fast^*}\; \!
\frakB_{\ol\calE,\ol\calR}(U,w; \Xi,\zeta)
 =  \frakB_{E,\calR_\slow}(U,\Xi)  + \scrB_\red(U,\Xi) 
\\
  \label{eq:Case2.scrLeff}
\text{with } &\scrB_\red(U,\Xi): =
 \sup_{w\in X_{\fast_{}} \atop P_\fast w=P_\slow U} 
  \inf_{\zeta \in X^*_{\fast_{}} \atop \PPC_\fast \zeta = \PPC_\slow \Xi}\! \!
    \frakB_{e,\calR_\fast}(w,\zeta) .
\end{align}

Thus, we see that $\scrB_\red$ is exactly obtained as in Section
\ref{su:ReducLagrDual}. Hence, we know that $\scrB_\red$ has a BER 
structure if for all $\eta\in Y$ the CSPP \eqref{eq:CSPP} for
$\frakB_{e,\calR_\fast}$ with constraint $\PPC \rmD e(w)=-\eta$ has a
null-saddle. In that case we have the BER structure $(\sfE_Y,\sfR_Y)$
such that
\[
 \scrB_\red (U,\Xi) = \frakB_{E,\calR_\red} (U,\Xi) \quad 
\text{with } \calR_\red^* (U,\Xi) = \sfR^*_Y\big( P_\slow U,
\PPC_\slow \Xi\big) .  
\]
We see that $\scrB_\red$ depends on $(U,\Xi)$ only through the
port values $\big( P_\slow U,  \PPC_\slow \Xi\big)\in Y\ti Y^*$. 
Returning to $\scrB_\eff = \frakB_{E,\calR_\slow} + \scrB_\red$ we obtain 
\[
\scrB_\eff = \frakB_{E, \calR_\eff} \quad \text{with } \calR^*_\eff(U,\Xi)= 
\calR_\slow(U,\Xi) + \sfR_Y\big( P_\slow U,  \PPC_\slow \Xi\big) . 
\]
Moreover, we see that $(X_\slow,E,\calR_\eff)$ is the EDP limit of
$(X,\calE_\eps,\calR_\eps)$ and the effective gradient-flow equation reads
\[
\dot U=\rmD_\Xi \calR^*_\eff \big(U,-\rmD E(U)\big) = \rmD_\Xi \calR^*_\slow
\big(U,-\rmD E(U)\big) + \PPO _\slow \rmD_\eta \sfR_Y\big( P_\slow U,
{-}\PPC_\slow \rmD E(U) \big) ,
\]  
which clearly shows that the non-equilibrium flux is given by 
\[
 \PPO_\slow V \quad \text{with } V= \rmD_\eta \sfR_Y\big( P_\slow U,
 \PPC_\slow \Xi \big) \in Y.
\]

\section{EDP-convergence for two ODE examples}
\label{se:ODEExamples}

We first treat the linear case as given in \eqref{eq:calLquadr} and with a
suitable scaling in $\eps>0$. Secondly, we consider a nonlinear reaction
systems with four species and two binary reactions $A+B \rightleftharpoons D$
and $A+ D  \rightleftharpoons C$ and show that the limiting system gives the
single ternary reaction $ 2A+B  \rightleftharpoons C$.

\subsection{Simple quadratic energy and dissipation}
\label{su:QuadrEnerDiss}

On the Hilbert space $X=X_\slow\ti X_\fast$ with $u=(U,w)$ we consider the family
$(X,\calE_\eps,\calR_\eps)$ of GSs given by $ \calE_\eps(U,w) = E(U)+ \eps
e(w)$ with 
\[
 E(U)=\frac12\langle
\bbA_\rms U-\mu_\rms,U\rangle_{X_\slow} \ \text{ and } \
  e(w)=\frac12\langle
\bbA_\rmf w-\mu_\rmf,w \rangle_{X_\fast}
\]
and 
\[
\calR_\eps(\Xi,\xi)=\frac12\Big\langle \binom{\Xi}{\frac1\eps\xi},
\binom{\bbK_{\rms\rms}\ \ \bbK_{\rms\rmf}}
{\bbK_{\rmf\rms}\ \ \bbK_{\rmf\rmf}} \binom{\Xi}{\frac1\eps\xi} \Big\rangle
= \ol\calR^*\big(\Xi,\frac1\eps \xi\big).
\]
Hence, we have the situation treated in Section \ref{su:Case1Product}.

The linear gradient-flow equation  takes the form 
\[
\binom{\dot U_\eps}{\eps \dot w_\eps} = - \binom{\bbK_{\rms\rms}\ \ \bbK_{\rms\rmf}}
{\bbK_{\rmf\rms}\ \ \bbK_{\rmf\rmf}} \binom{\bbA_\rms
  U_\eps-\mu_\rms}{\bbA_\rmf w_\eps-\mu_\rmf} .
\]
With the port mappings $P(U,w)=U \in Y:=X_\slow$ and
$ \PPC (\Xi,\zeta) \to \Xi \in X_\fast $ we obtain the determining equation
\eqref{eq:NESSeqn} for the NESS
\[
\binom{0}{0} = - \binom{\bbK_{\rms\rms}\ \ \bbK_{\rms\rmf}}
{\bbK_{\rmf\rms}\ \ \bbK_{\rmf\rmf}} \binom{ \rmD E(U)}{\rmD e(w)} +
\binom{V}{0}, \quad \rmD E(U) = -\Xi \in X_\slow , \quad V\in X_\slow. 
\]
As $\Xi$ is given, and the upper equation is always true for a suitable $V$, we
find the NESS 
\[
 \bbA_\rmf w {-} \mu_\rmf = \rmD e(w) = \bbK_{\rmf\rmf}^{-1} \bbK_{\rmf\rms} \Xi.
\]
The resulting port mapping $\mfP:X^*_\slow \to X_\slow;\ \Xi\mapsto V$ takes
the explicit form 
\[
V= \mfP \Xi =  \bbK_\eff \Xi \quad \text{with }  \bbK_\eff
= \bbK_{\rms\rms}- \bbK_{\rms\rmf}\bbK_{\rmf\rmf}^{-1}\bbK_{\rmf\rms}.
\]
In particular, $\mfP$ is independent of the energy $\calE$, as
predicted by Proposition \ref{pr:PortRelaSIDiss}.

We also want to show that $\mfP=\rmD\calR^*_\eff$ can be obtained by the
saddle-point reduction of the B-function
\[
\frakB_{\ol\calE,\ol\calR}(U,w;\Xi,\zeta) = \ol\calR^*(\Xi,\zeta) -
\ol\calR^* \big(\mu_\rms{-} \bbA_\rms U, \mu_\rmf{-} \bbA_\rmf w\big).
\]
Assuming that $\bbK>0$ and $\bbA_\rms>0$, a simple calculation gives 
\begin{align*}
\scrB_\red(U,\Xi)& = \! \sup_{w\in X_\fast} \inf_{\zeta \in X_\fast^*} \!\!
\frakB_{\ol\calE,\ol\calR}(U,w;\Xi,\zeta) = \! \inf_{\zeta \in
  X_\fast^*} \!\! \ol\calR^*(\Xi,\zeta) - \!\! \inf_{ w\in X_\fast} \!\! 
\ol\calR^* \big(\mu_\rms{-} \bbA_\rms U, \mu_\rmf{-} \bbA_\rmf w\big) 
\\
& = \frac12 \langle \Xi, \bbK_\eff \Xi\rangle - \frac12 \langle \mu_\rms{-}
\bbA_\rms U, \bbK_\eff (\mu_\rms{-} \bbA_\rms U)\rangle = \frakB_{E,\calR_\eff}(U,\Xi)
\end{align*}
with $\calR_\eff(\Xi)= \frac12 \langle \Xi, \bbK_\eff \Xi\rangle$.

\subsection{Two binary reaction generate one ternary reaction}
\label{suBinaryTernaryReact} 

We consider four chemical species $A$, $B$, $C$, and $D$ with associated
concentrations $a,b,c,d\in {[0,\infty[}$. They undergo the two binary reversible
reaction pairs $A+B \rightleftharpoons D$ and $A+D \rightleftharpoons C$
according to the mass action law.  We assume that species $D$ is very unstable
and either react fast with an $A$ to create $C$ or decay fast into $A$ and
$B$. In particular, the equilibrium concentrations for $D$ will be
$d_\eps := \eps w_*$, while the equilibrium densities $a_*,b_*, c_*$ are
positive and independent of $\eps$, see Figure \ref{fig:RRE}. 
\begin{figure}
\centering 
\begin{tikzpicture}
\tikzstyle{species} = [circle, fill=gray!30];

\draw[thick, fill=gray!10, rounded corners] (-0.6,-1.3) rectangle (5.5,2.4);
\node[species] at (0,1.5) {A} ;
\node[species] at (0,0) {B};
\node[species] at (2.5,0.75) {D};
\node at (0,0.75) {$\bm+$}; 
\draw [very thick, ->] (0.4,0.88)-- 
      node[above, pos=0.5]{{\scriptsize slow}}  (2,0.88);
\draw [ultra thick, ->] (2.,0.62)-- 
      node[below, pos=0.5]{{\scriptsize fast}}  (0.4,0.62);
\node[species] at (2.5,-0.75) {A} ;
\node at (2.5,0.0) {$\bm+$}; 
\draw [very thick, ->]  (4.5,-0.13) -- 
      node[below, pos=0.5]{{\scriptsize slow}} (2.9,-0.13);
\draw [ultra thick, ->]  (2.9,0.13) -- 
      node[above, pos=0.5]{{\scriptsize fast}} (4.5,0.13);
\node[species] at (5.0,0) {C};
\node at (2.5,2) {slow-fast system};

\draw[ultra thick,->] (6,0.65)-- node[above, pos=0.5] {$\eps \to 0$} (8,0.65);

\begin{scope}[shift={(10,-0.4)}]
\draw[thick, fill=gray!10, rounded corners] (-1.5,-.55) rectangle (3.1,2.6);
\node[species] at (0,1.5){A} ; 
\node[species] at (0,0){A} ;
 \node[species] at (-0.9,0.75){B} ;
\node at (-0.2 ,0.75) {$\bm+$}; 
\draw [very thick, ->] (0.4,0.88)-- 
      node[above, pos=0.5]{{\scriptsize slow}}  (2,0.88);
\draw [very thick, ->] (2.,0.62)-- 
      node[below, pos=0.5]{{\scriptsize slow}}  (0.4,0.62);
\node[species] at (2.5,0.75) {C};
\node at (1,2.2) {effective system};
\end{scope}
\end{tikzpicture}

\caption{The slow-fast reaction-rate equation \eqref{eq:RRE.bin.bin} has two
  binary reaction pairs with one fast and one slow reaction. The effective
  system \eqref{eq:RRE.ter} has one slow ternary reaction pair.}
\label{fig:RRE}
\end{figure}

The associated reaction rate equation is the ODE system 
\begin{equation}
  \label{eq:RRE.bin.bin}
  \bma{c} \dot a\\ \dot b \\ \dot c \\ \dot d\ema = 
\kappa_1 \Big( \frac{d}{d_\eps} - \frac{ab}{a_*b_*}\Big) 
  \bma{c} 1\\ 1\\ 0 \\ \!\!-1\!\! \ema  
+ \kappa_2 \Big( \frac{c}{c_*} - \frac{ad}{a_*d_\eps}\Big)
   \bma{c} 1\\ 0\\ \!\!-1\!\! \\ 1 \ema , 
\end{equation}
where $\kappa_1$ and $\kappa_2$ are positive reaction coefficients that may
depend on $a,b,c,d$, but make them constant for simplicity. 

As above one may replace $d$ by $\eps w$ and such that the right-hand side
becomes independent of $\eps$. Dropping the term $\eps w$ on the left-hand side
leads to the algebraic-differential system 
\begin{equation}
  \label{eq:RRE.bin2.eps0}
  \bma{c} \dot a\\ \dot b \\ \dot c \\ 0 \ema = 
\kappa_1 \Big( \frac{w}{w_*} - \frac{ab}{a_*b_*}\Big) 
  \bma{c} 1\\ 1\\ 0 \\ \!\!-1\!\! \ema  
+ \kappa_2 \Big( \frac{c}{c_*} - \frac{aw}{a_*w_*}\Big)
   \bma{c} 1\\ 0\\ \!\!-1\!\! \\ 1 \ema , 
\end{equation}
Solving the last equation for $w$ and inserting the result into the first three
equations leads to the reduced ODE 
\begin{equation}
  \label{eq:RRE.ter}
  \bma{c} \dot a\\ \dot b \\ \dot c \ema =  \kappa_\eff(a) \Big( \frac{c}{c_*} -
  \frac{a^2b}{a^2_*b_*}\Big) \bma{c} 2\\ 1\\  \!\!-1\!\! \ema \quad \text{with } 
  \kappa_\eff(a) := \frac{\kappa_1\kappa_2 a_*}{\kappa_1 a_* {+} \kappa_2 a}
\end{equation}
which is the reaction-rate equation for the ternary reaction $2A+B
\rightleftharpoons C$ with an effective reaction coefficient $\kappa_\eff(a)\in
{]0,\kappa_2[}$. 

The original system has the entropic cosh-gradient structure as derived in
\cite{MPPR17NETP} and further studied in \cite{MieSte20CGED, MiPeSt21EDPC}.
In our specific case, the reaction-rate equation \eqref{eq:RRE.bin.bin} 
is the gradient-flow equation for the GS $(\R^4,\calE_\eps,\calR_\eps)$ given by
(where $u=(a,b,c,d)$) 
\begin{equation}
  \label{eq:EBz.R*cosh}
  \begin{aligned}
&\calE_\eps(u)=\LB(a/a_*)a_*+ \LB(b/b_*)b_*+\LB(c/c_*)c_*+\LB(d/d_\eps) d_\eps 
\ \text{ and } \\
& 
\calR^*_\eps(u;\xi) = \kappa_1\big(\frac{abd}{a_*b_*d_\eps}\big)^{1/2}
\sfC^*\big(\xi_1{+}\xi_2{-}\xi_4\big) + \kappa_2\big(\frac{acd}{a_*c_*d_\eps}\big)^{1/2}
\sfC^*\big(\xi_1{-}\xi_3{+}\xi_4\big),
\end{aligned}
\end{equation}
where $\LB(z)=z \log z -z+1$ is the Boltzmann function and $\sfC^*(\zeta)=4
\cosh(\zeta/2)-4$. 

Doing our standard scaling for the slow and fast variables gives 
\[
u=(U,\eps w), \quad U=(a,b,c)\in X_\slow, \quad \calE_\eps(u) =
E(U) + \eps e(w) \text{ with } e(w)=\LB(w/w_*)w_*.
\]
Moreover, with $\Xi=(\xi_1,\xi_2,\xi_3) \in X_\slow^*$ we have
$\wt\calR_\eps(U,\eps w;\Xi,\mu) = \ol\calR\big(U,w;\Xi,\frac1\eps \mu\big) $
with
\[
\ol\calR^*\big(U,w;\Xi,\zeta\big) 
 = \kappa_1\big(\frac{abw}{a_*b_*w_*}\big)^{1/2}
\sfC^*\big(\xi_1{+}\xi_2{-}\zeta\big) + \kappa_2\big(\frac{acw}{a_*c_*w_*}\big)^{1/2}
\sfC^*\big(\xi_1{-}\xi_3{+}\zeta \big),
\]

Thus, we can apply the theory of Section \ref{su:Case1Product} and define
$\scrB_\red $ as in \eqref{eq:def.calLslow}, namely
\[
\scrB_\red(U,\Xi):= \sup_{w>0} \inf_{\zeta\in \R}
\frakB_{\ol\calE,\ol\calR}(U,w;\Xi,\zeta). 
\]
The sup-inf can be calculated explicitly as is explained in
\cite[Sec.\,3.3.2]{LMPR17MOGG}. Indeed using the formula
\[
\inf_{\zeta\in \R} \big(g\sfC^*(\zeta)+ h\sfC^*(\rho-\zeta)\big) =
4W(g,h,\rho) -4(g{+}h) \text{ with } W(g,h,\rho)= \big((g{+}h)^2 + \frac{gh}2
\sfC^*(\rho)\big)^{1/2} ,
\]
where $\rho=2\xi_1{+}\xi_2{-}\xi_3$,  $g=
\kappa_1\big(\frac{abw}{a_*b_*w_*}\big)^{1/2}$, and 
$h =\kappa_2\big(\frac{acw}{a_*c_*w_*}\big)^{1/2} $, a lengthy calculation yields
\[
\scrB_\red(U,\Xi):= \sup_{w>0}\Big( 4W(g,h,\rho)
 -2 \kappa_1 \big(\frac{ab}{a_*b_*}    {+}\frac w{w_*} \big)  
-2 \kappa_2\big( \frac c{c_*} + \frac{aw}{a_*w_*} \big) \Big)
\] 
Noting that $g$ and $h$ are proportional to $\sqrt w$, we see that also
$W(g,h,\rho)$ is exactly proportional to $\sqrt w$. Hence, the maximum with
respect to $w$ can be determined and another lengthy calculation gives the
explicit expression
\[
\scrB_\red(U,\Xi) = \kappa_\eff(a)
\big(\frac{a^2\,b \,c}{a_*^2b_*c_*}\big)^{1/2} \,\sfC^*\big(
2\xi_1{+}\xi_2{-}\xi_3\big) - \kappa_\eff(a)\,2\, 
\Big( \big(\frac{a^2\,b}{a_*^2b_*} \big)^{1/2}- \big(\frac
c{c_*}\big)^{1/2} \Big)^2   
\]
with $\kappa_\eff(a)$ from \eqref{eq:RRE.ter}. 
Hence, we have the BER structure
$\scrB_\red(U,\Xi)= \calR_\eff^*(U,\Xi)-\calR^*_\eff (U,-\rmD E(U))$ with
$\calR_\eff^*(U,\Xi) = \kappa_\eff(a)
\big(\frac{a^2\,b \,c}{a_*^2b_*c_*}\big)^{1/2} \,\sfC^*\big(
2\xi_1{+}\xi_2{-}\xi_3\big)$. 

It seems that the above theory can be generalized to an arbitrary number of
species with a density vector $\bfc=(c_1,...,c_{i_*}) \in \R^{i_*}$ and an
arbitrary number $ r_*$ of reactions following the mass-action law, as long as
we have the detailed-balance condition, i.e.\ there exists a positive steady
state $\bfc^*_\eps= (c_1^*,...,c_{j_*},\eps w^*_{j+1},...,\eps w^*_{i_*})$. If
this is so, then the interesting question is how the reaction coefficients of
the limiting system depend on the reaction coefficients of the original
system. Note that even in our simple case, we can start with constant
coefficients $\kappa_1$ and $\kappa_2$ but then find $\kappa_\eff(a)$ which 
depends on the state.

In particular, we want to highlight that the effective system has again the
expected entropic cosh-gradient structure for the ternary reaction
$2A+B\rightleftharpoons C$. We emphasize that this is not automatic, because in
\cite[Sec.\,4.3]{MiPeSt21EDPC} an example of a reaction-rate equation is
studied where the EDP-limit of the entropy cosh-gradient structure leads to an
effective GS $(\R^4,\bfE,\bfR)$ where $\bfE$ is no longer a
Boltzmann entropy and the reaction does no longer follow the mass-action law.

\section{Linear diffusion through a membrane}
\label{se:LinearDiffMembrane}

The example in this section is well studied from the context of PDEs and
singular limits. We are looking at a diffusion system of $i_*$ mass densities
$\bfrho = (\rho_1,...,\rho_{i_*})$ that diffuse along an interval on the real
line, where in the small interval ${]{-}\eps,\eps[}$ representing a membrane 
the mobility is also of
order $\eps$, whereas it is of order $1$ outside the membrane. 

\subsection{The PDE model and its quadratic gradient structure}
\label{su:DiffMembrPDE}

We consider the intervals $\Omega_\eps={]{-}1{-}\eps, 1{+}\eps[}$ and define
the piecewise affine maps $\psi_\eps$ and $\phi_\eps$ between $\Omega_\eps$ and
$\Omega:=\Omega_1={]{-}2,2[}$:
\begin{equation}
  \label{eq:psieps.phieps}
  \psi_\eps(x)=\begin{cases}  x{+}\eps{-}1&\text{for } x\geq 1, \\ \eps\, x&
  \text{for } |x|\leq 1, \\ x{-}\eps{+}1& \text{for } x\leq -1; 
\end{cases}
\quad \text{and} \quad 
\phi_\eps(y)=\begin{cases}  y{-}\eps{+}1&\text{for } y\geq \eps, \\ y/\eps &
  \text{for } |y|\leq \eps, \\ y{+}\eps{-}1& \text{for } y\leq -\eps, 
\end{cases}
\end{equation}
see Figure \ref{fig:psi.phi.Keps}. 
The original diffusion problem is defined on $\Omega_\eps$ and we assume that
the mobility is given in the form
\begin{equation}
  \label{eq:def.PC0}
  K_\eps (y) = \frac1{\phi'_\eps(y)} \,\ol K \big(\phi_\eps(y)\big) \quad
\text{with }
\ol K\in \rmP\rmC^0 \big( [-2,-1] \,\ol\cup \, [-1,1] \,\ol\cup\, [1,2];\R^{i_*\ti
  i_*}_\mafo{sym} \big),
\end{equation}
where the crucial point is that the prefactor $1/\phi'_\eps(y)$ provides
the prefactor $\eps$ inside the membrane region ${]{-}\eps,\eps[}$. The
notation $\rmP\rmC^0$ with ``$\ol\cup$'' means that $\ol K$ is piecewise
continuous and has continuous extensions on the three closed intervals
$[-2,-1]$, $[-1,1]$, and $[1,2]$, such that the one-sided limits
\[
K_\pm := \ol K\big(\pm (1{+}0)\big) = \lim_{\delta\to 0^+} \ol
K\big(\pm(1{+}\delta) \big) \ \text{ and } \ 
k_\pm := \ol K\big(\pm (1{-}0)\big) = \lim_{\delta\to 0^+} \ol K\big(\pm(1{-}\delta) \big)
\]
exist, but may be different.  Moreover, we assume that $\ol K$ is positive definite,
i.e.\ there exists $ \kappa >0$ such that $a\vdot \ol K(x) a\geq \kappa |a|^2$
for all $x\in [-2,2]$ and $a\in \R^{i_*}$. Hence, $y\mapsto K_\eps(y)$ is
discontinuous at $y=\pm \eps$, because it jumps by a factor of $\eps$.
\begin{figure}
\newcommand{\EPS}{0.1}
\centering
\begin{tikzpicture} % eps = 0.1 

\draw[gray!10, fill=gray!10] (-3,-2.8) rectangle (11,2.7);
\draw[gray,thick,->] (-2.25,0) -- (2.25,0) node[above]{$x$};
{\footnotesize\draw[gray] (-2,0.1)--(-2,-0.1) node[below]{$-2$};
\draw[gray] (-1,0.1)--(-1,-0.1) node[below]{$-1$};
\draw[gray] (1,0.1)--(1,-0.1) node[below]{$1$};
\draw[gray] (2,0.1)--(2,-0.1) node[below]{$2$};}
\draw[gray,thick,->] (0,-1.25)--(0,1.25) node[left]{$y$};
\draw[ultra thick] (-2,-\EPS-1) -- (-1,-\EPS) -- (1,\EPS)
    -- node[pos =0.5, above left] {$\psi_\eps(x)\!\!$} (2,1+\EPS);
\draw[thick, |-|] (-2,-1.3)-- node[pos=0.5,
below]{$\Omega_\slow$} (-1,-1.3);
\draw[thick, |-|] (-1,-2)-- node[pos=0.5, below]{$\Omega_\fast=[-1, 1]$} (1,-2);
\draw[thick, |-|] (1,-1.3)-- node[pos=0.5, below]{$\Omega_\slow$} (2,-1.3);

\begin{scope}[shift={(4.5,-0.0)}]
\draw[gray,thick,->] (-1.25,0) -- (1.5,0) node[above]{$y$};
{\footnotesize\draw[gray] (-0.1,-2)--(0.1,-2) node[right]{$-2$};
\draw[gray] (-0.1,-1)--(0.1,-1) node[right]{$-1$};
\draw[gray] (0.1,1)--(-0.1,1) node[left]{$1$};
\draw[gray] (0.1,2)--(-0.1,2) node[left]{$2$};}
\draw[gray,thick,->] (0,-2)--(0,2.25) node[right]{$x$};
\draw[ultra thick] (-\EPS-1,-2) -- (-\EPS,-1) -- (\EPS,1)
    -- node[pos =0.5, below right] {$\!\!\phi_\eps(y)$} (1+\EPS,2);
\end{scope}

\begin{scope}[shift={(8.5,-0.0)}]
\draw[gray,thick,->] (-1.25,0) -- (1.35,0) node[above]{$y$};
{\footnotesize
\draw[gray] (-1,0.1)--(-1,-0.1) node[below]{$-1$};
\draw[gray] (1,0.1)--(1,-0.1) node[below]{$1$};}
\draw[gray,thick,->] (0,-0.25)--(0,2.0) node[left]{$K_\eps(y)$};
\draw[ultra thick, domain=-1-\EPS:-\EPS, samples =30, smooth, variable=\x]
                         plot ({\x},{1+0.6*cos(\x r)});
\draw[ultra thick, domain=-\EPS:\EPS, samples =10, smooth, variable=\x]
                         plot ({\x},{\EPS+\EPS*sin(\x r)});
\draw[ultra thick, domain=\EPS:1+\EPS, samples =30, smooth, variable=\x]
                         plot ({\x},{1.1-0.6*sin(2*\x r)});

\draw[thick, |-|] (-1-\EPS,-1)-- node[pos=0.5,
below]{$\Omega_\eps=[-1^{}{-}\eps, 1{+}\eps]$} (1+\EPS,-1);
\end{scope}
\end{tikzpicture}

\caption{Left and middle: the diffeomorphisms $\psi_\eps$ and $\phi_\eps$ map
  $\Omega_\fast\cup \Omega_\slow= [-2,2]$ to $\Omega_\eps$. Right: the positive
  mobility $K_\eps$ is order $\eps$ in the membrane ${]{-}\eps,\eps[}$ and
  order $1$ otherwise.}
\label{fig:psi.phi.Keps}
\end{figure}

\begin{figure} 
\centering
\begin{tabular}{cc}
\includegraphics[width=0.4\textwidth, trim=0 30 0 30, clip=true]{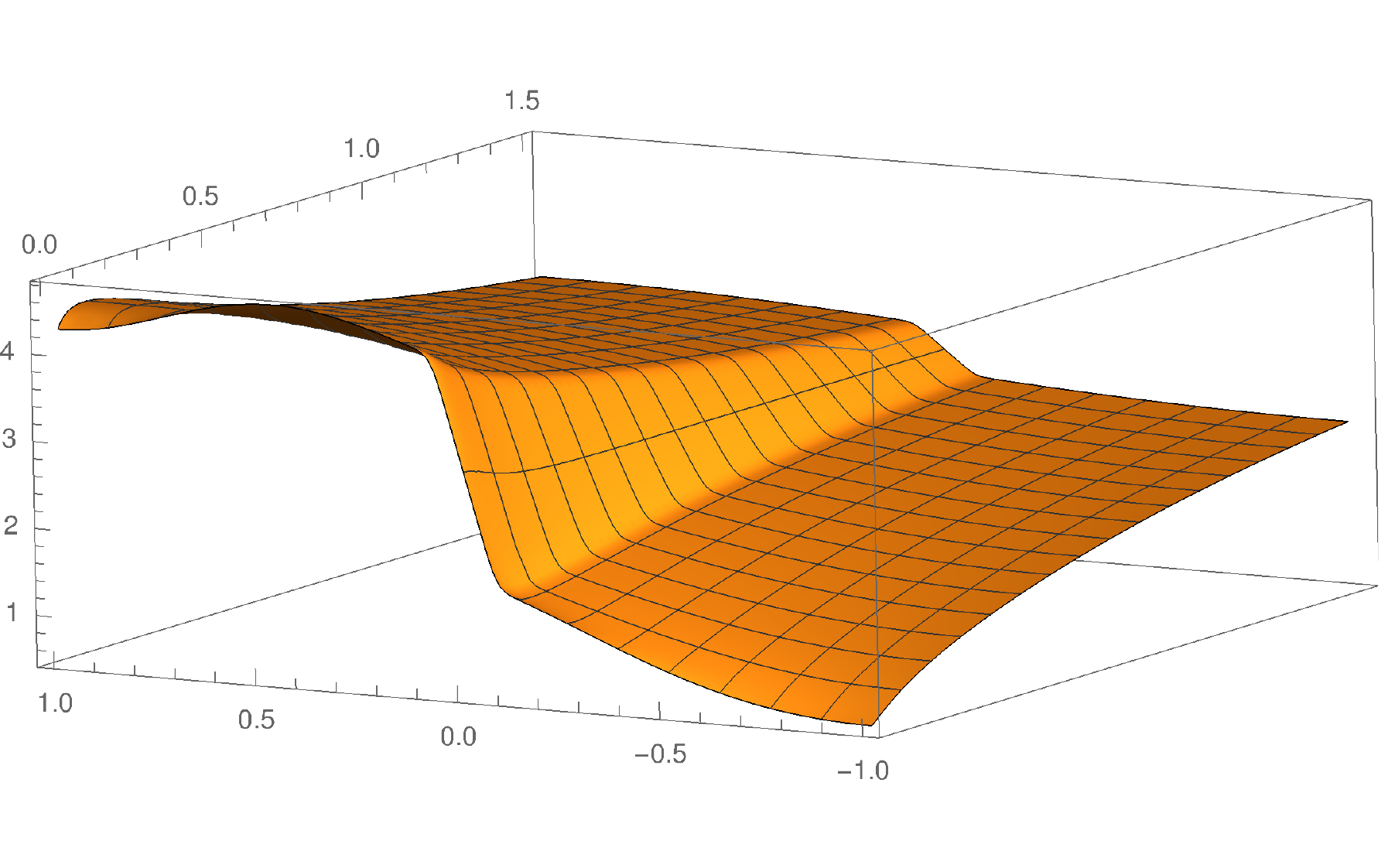} &
\includegraphics[width=0.4\textwidth]{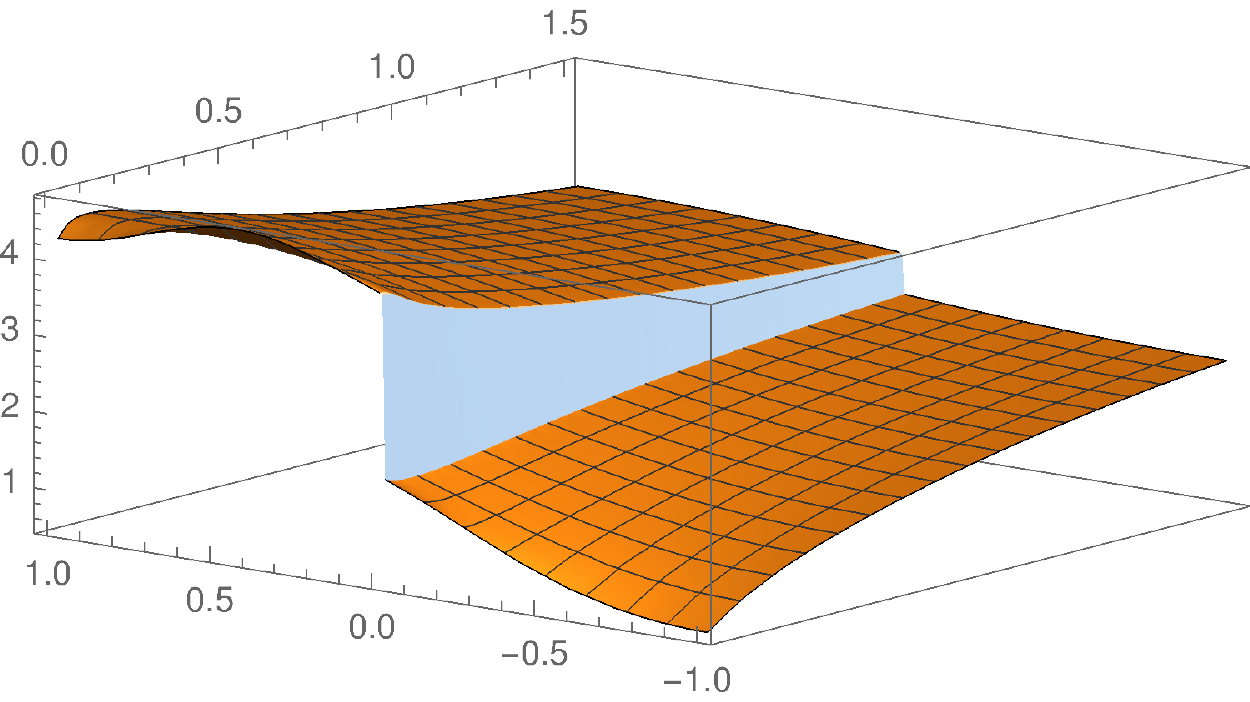} 
\\
slow-fast diffusion with $\eps=0.1$  & effective diffusion/transmission
for$\eps\to0$
\end{tabular}
\caption{Left: the solution of \eqref{eq:DiffMembEps} for $\eps=0.1$ shows a
  steep slope in the membrane ${]{-}\eps,\eps[}$. Right: the
  solution of the effective transmission problem \eqref{eq:DiffMembEff} jumps at $y=0$.}
\label{fig:Sol.eps.0}
\end{figure}

We define a second positive definite function
$\ol A\in \rmP\rmC^0 \big( [-2,-1]\,\ol\cup \,[-1,1]\,\ol\cup\,[1,2];\R^{i_*\ti
  i_*}_\mafo{sym} \big)$ which determines the energy functional
\[
\wt\calE_\eps( \bfrho) :=\int_{\Omega_\eps} \frac12 \bfrho(y) \vdot \ol
A(\phi_\eps(y)) \bfrho(y) \dd y \text{ \ on the space } X_\eps
 =\rmL^2(\Omega_\eps; \R^{i_*}).
\]
Moreover, we define the dual dissipation potential $\calR_\eps$ via
\[
\wt\calR_\eps(\bfmu) = \int_\Omega \frac12 \pl_y\bfmu(y) \cdot K_\eps(y)
\pl_y\bfmu(y) \dd y.
\]
The gradient-flow equation for the GS $(X_\eps,\wt\calE_\eps,\wt\calR_\eps)$ is the linear
parabolic system  
\begin{equation}
  \label{eq:DiffMembEps}
  \dot\bfrho= \pl_y\Big( K_\eps(y) \,\pl_y\big( A_\eps(y) \bfrho(t,y)\big) \Big) \
\text{ for } t>0,\ y\in \Omega_\eps, \quad \pl_y\big(
A_\eps(y)\bfrho(t,y)\big)\big|_{y=\pm(1{+}\eps)} =0 .
\end{equation}
Note that $\bfM(t)=\int_{\Omega_\eps}\bfrho(t,y) \dd y \in \R^{i_*}$ is
independent of $t$ because of the divergence form and the no-flux boundary
conditions.  A typical solution with large gradients for $y \in
{]-\eps,\eps[}$ is depicted in Figure \ref{fig:Sol.eps.0}. 

To study the limit $\eps\to 0$ it is advantageous to transform the PDE to the
fixed interval $\Omega$ via $\phi_\eps(\Omega_\eps)=\Omega$. 
For $x \in \Omega$ we set 
\[
\bfu(t,x)=\frac1{\psi'_\eps(x)}\,\bfrho(t, \psi_\eps(x)) \quad 
 \text{and} \quad \calE_\eps(\bfu) 
= \wt\calE_\eps\big(\tfrac1{\phi'_\eps}\bfu{\circ}\phi_\eps\big) 
 = \int_\Omega \frac12 \bfu\vdot \ol
    A \bfu\, \psi'_\eps\, \dd x.
\]
The transformed dissipation potential takes the form 
\[
\calR_\eps(\bfxi) =\int_\Omega \frac12 \pl_x\big(\frac1{\psi'_\eps}\bfxi\big)
\vdot \ol K\, \pl_x\big( \frac1{\psi'_\eps}\bfxi\big)  \,\dd x,
\]
where we used the scaling $\phi'_\eps(y)K_\eps(y)=\ol K(x)$ to cancel the
powers of $\psi'_\eps$. 

The transformed  linear diffusion equation reads
\[
\psi'_\eps(x) \,\dot\bfu (t,x)= \pl_x\Big( \ol K(x)\,\pl_x\big( \ol A(x)
\bfu(t,x)\big) \Big) , \quad \pl_x\big( \ol A(x)
\bfu(t,x)\big)\big|_{x=\pm 2} = \bm0. 
\]
Of course, in the above development we have anticipated the scalings in such a
way that in the last equation $\eps$ only occurs once, namely in the prefactor
$\psi'_\eps$ with $\psi'_\eps(x)=\eps$ for $|x|<1$ and $\psi'_\eps(x) = 1$
for $1<|x|<2$. Thus, we are exactly in the situation of a slow-fast gradient
system as studied in Section \ref{se:ReductGSNESS}. 

We make the splitting and the corresponding port mappings explicit. We are in
``Case 2'' where the product space $X=X_\slow\ti X_\fast$ needs a factorization
along the boundary of the membrane, now placed at $x= \pm 1$. We set (see
Figure \ref{fig:psi.phi.Keps}) 
\[
\Omega_\fast = {[{-}1,1]}, \quad \Omega_\slow ={]{-}2,-1]}\cup {[1,2[}, \quad 
X_\fast = \rmL^2(\Omega_\fast;\R^{i_*}), \quad X_\slow =
\rmL^2(\Omega_\slow;\R^{i_*})
\]
and introduce the variable $ U = u|_{\Omega_\slow}\in X_\slow $ and $ w = 
u|_{\Omega_\fast} \in X_\fast $. With this we find the transformed energy 
\[
\calE_\eps(U,w) = E(U)+ \eps \,e(w) \quad \text{with } E(U)=\int_{\Omega_\slow}
  \!\! \frac12 U\cdot \ol A U \dd x \text{ and } 
e(w) = \int_{\Omega_\fast} \!\! \frac12 w\cdot \ol A w \dd x. 
\]
If we similarly write $\bfxi =(\Xi,\zeta)$ with $\Xi=\bfxi|_{\Omega_\slow} \in
X^*_\slow$ and $\zeta= \frac1\eps \bfxi|_{\Omega_\fast}\in X^*_\fast$ we obtain 
\begin{equation}
\label{eq:DiffMembCompat}
\begin{aligned}
&\ol\calR^*(\Xi,\zeta) = \calR^*_\slow(X) + \calR^*_\fast(\zeta) +
\bfdelta_{\{0\}}\big( \PPC_\slow \Xi{-} \PPC_\fast \zeta\big) 
\\
&\text{where }
\calR^*_\slow(\Xi)=\int_{\Omega_\slow}\!\!\frac12 \pl_x \Xi \vdot \ol K
\pl_x \Xi \dd x \text{ and } 
\calR^*_\fast(\zeta)=\int_{\Omega_\fast}\!\!\frac12 \pl_x \zeta \vdot \ol K
\pl_x \zeta \dd x. 
\end{aligned}
\end{equation}
Here the compatibility condition $ \PPC_\slow \Xi = \PPC_\fast \zeta$ is
crucial. We define $Y=\R^{i_*}\ti \R^{i_*}$ and the port mappings (where
$f(x^+)$ and $f(x^-)$ denote the limit from the right and left, respectively)
\[
P_\slow:X_\slow \to Y; U\mapsto (U(-1^-),U(1^+)) \  \text{ and } \ 
P_\fast:X_\fast \to Y; w\mapsto (w(-1^+),w(1^-)), 
\]
and similarly $\PPC_\slow :X^*_\slow \to Y^*$ and $\PPC_\fast:X^*_\fast \to Y^*$. 

The limiting equation for $\eps=0$ takes the form 
\begin{align}
\nonumber
\dot U&= \pl_x\big( \ol K\,\pl_x(\ol A U) \big) \ \text{ for }x\in {]1,2[}, 
   &&\pl_x(\ol A U)|_{ x=2} =0,
\\\nonumber
0&=  \pl_x\big( \ol K\,\pl_x(\ol A w) \big) \ \text{ for }x\in {]{-}1,1[}, 
&&\hspace*{-1em}\left\{\begin{aligned} &\ \ \,U(1^+)=w(1^-), \quad \ \ 
   \pl_x(\ol A U)|_{x=1^+}\!= \pl_x(\ol A w)|_{x=1^-}, 
  \\ 
  &U({-}1^-)=w({-}1^+), \  \pl_x(\ol A U)|_{x={-}1^-} \!= \pl_x(\ol A w)|_{x={-}1^+},
  \end{aligned} \right.
  \\
\label{eq:DiffMembEff}
\dot U&= \pl_x\big( \ol K\,\pl_x(\ol A U) \big) \ \text{ for }x\in {]{-}2,{-}1[}, 
   &&  \pl_x(\ol A U)|_{ x=-2} =0.
\end{align}
The static equation on $\Omega_\fast={[{-}1,1]}$ can be solved explicitly 
via $\ol K(x)\pl_x\big(\ol A(x)w(x)\big)=\mafo{const.}$, and we
obtain the corresponding \emph{transmission conditions}
\[
\ol K\pl_x(\ol A U)\big|_{x=\pm 1} = \bbH_K \big( A(1)U(1)- A(-1)U(-1)\big), \
\text{ where }\bbH_K=\Big( \int_{-1}^1 \ol K(x)^{-1}\dd x \Big)^{-1}.
\]

To understand the thermodynamical origin of the effective transmission
conditions, we use 
EDP-convergence via reduced B-functions as described in Section
\ref{su:Case2Product}. For this, we  construct
\[
\scrB_\eff (U,\Xi)
= \sup \inf \frakB_{\ol\calE,\ol\calR} (U,w;\Xi,\zeta) =
\frakB_{E,\calR_\slow}(U,\Xi) + \scrB_\red (U,\Xi) 
\]
We are in the case where $\ol\calR$ is independent of the state, such that
$\calR_\red$ has the form 
\[
\calR^*_\red(\Xi)= \sfR_Y^*(\PPC_\slow \Xi) \quad \text{with } 
\sfR_Y^*(\eta):=\inf_{\Xi:\; \PPC_\slow \Xi=\eta} \calR^*_\fast(\Xi),
\]
see Proposition \ref{pr:PortRelaSIDiss}.  A direct calculation shows
that $\mfP$ is given in terms of 
\[
\calR^*_\red(\Xi)=\sfR^*_Y(\PPC_\slow \Xi) \quad \text{with }
\sfR_Y(\eta(-1),\eta(1)) = \frac12 \big(\eta(1){-}\eta(-1)\big) \vdot \bbH_K
\big(\eta(1){-}\eta(-1)\big) 
\]
which shows
$\mfP(\eta(1),\eta(-1)) =\big(\, \bbH_K (\eta(1){-}\eta(-1))\, , \,\bbH_K
(\eta(-1){-}\eta(1)) \, \big)$.  Indeed, $\calR_\red$ can easily be obtained by
minimizing $\calR_\fast(\zeta)$ over the constraints
$\PPC_\fast \zeta= \PPC_\slow \Xi$.

In summary, the effective GS $(X_\slow,E,\calR_\eff)$ is given by 
\[
\calR^*_\eff(\Xi)=\int_{-2}^{-1}\! \frac12 \Xi'\cdot \ol K \Xi' \dd x +
\frac12\big(\Xi(1) {-} \Xi(-1)\big)\cdot \bbH_K\big(\Xi(1) {-} \Xi(-1)\big) + 
 \int_{1}^{2}\! \frac12 \Xi'\cdot \ol K \Xi' \dd x. 
\]
We note that $\calR^*_\eff$ is independent of $\ol A$ (from the energy), which
is in contrast to the result using the Otto gradient structure. 

\subsection{EDP-convergence in the Otto gradient structure}
\label{su:DiffMembrOtt}

We reconsider the above linear equation, but now we strict to the
scalar case $i_*=1$, viz.\ $u(t,x)\in {[0,\infty[}\in \R^1$. The
linear equation can then be interpreted as a Fokker-Planck equation. Our aim is
to do the EDP-limit $\eps\to 0$ as in the previous subsection, but now
for the so-called Otto gradient structure, also called gradient-flow in the
Wasserstein space. The gradient system is the triple
$(\mafo{Prob}(\Omega_\eps), \calE^\rmB_\eps, \calR_\eps^\mafo{Otto})$, where
the function space is
\[
\mafo{Prob}(\Omega_\eps):=\bigset{ u \in \rmL^1(\Omega_\eps)}{ u\geq 0,\
  \int_{\Omega_\eps} u\dd y =1},
\]
the energy is Boltzmann's relative entropy
\[
\calE^\rmB_\eps(u)= \int_{\Omega_\eps} \LB\big( A_\eps(y)u(y)\big) \frac{\dd y
}{A_\eps(y)} ,
\]
and the dual dissipation functional reads
\[
\calR^\mafo{Otto\:*}_\eps (u,\xi) = \int_{\Omega_\eps} \frac{K_\eps(y)}2
|\pl_y\xi(y)|^2 u(y)\dd y = \frac12\big\langle \bbK_\eps^\mafo{Otto}(u) \xi,\xi\rangle,
\]
which is quadratic in $\xi$ and dependent on the state $u \in
\mafo{Prob}(\Omega_\eps)$. Here $\bbK_\eps^\mafo{Otto}(u)$ can be understood as
the self-adjoint nonnegative differential operator
\[
\bbK_\eps^\mafo{Otto}(u) \,\xi = - \pl_y \big( K_\eps u \pl_y\xi\,\big) \quad
\text{with } K_\eps u \pl_y\xi\big|_{y=\pm(1{+}\eps)}=0.
\]
The associated gradient-flow equation is the Fokker-Planck equation
\[
\dot u = \pl_y\big( K_\eps u \pl_y(A_\eps u)\big) = \pl_y\big( K_\eps(\pl_y u+
u V'_\eps)\big), 
\]
if we define the driving potential $V_\eps$ by $V_\eps(y) = \log A_\eps(y)$.
We refer to \cite{Otto96DDDE, JoKiOt97FEFP, Otto98DLPF, JoKiOt98VFFP,
  Otto01GDEE} for the first work treating the Fokker-Planck equation as an
gradient-flow equation with respect to this gradient structure.

We now want to do the EDP-limit in this gradient structure, where the new
feature is the dependence of $\calR^*_\eps$ on the state $u$. As a result the limit
gradient structure will be quite different. First it will depend in properties
of $\ol A$ which shows that $\calR^*_\slow$ cannot be calculated from
$\calR_\eps^*$ alone. Secondly, we will see that $\calR^*_\eff=\calR_\slow^* +
\calR^*_\red$ will no longer be quadratic in $\xi$, namely $\calR_\red^*$, which
is obtained from the NESS problem of the rescaled membrane, will have a
cosh-type behavior given through $\sfC^*$. 

We will not give the analysis in detail, as the result is well-established see
\cite[Sec.\,4]{LMPR17MOGG}, \cite[Sec.\,4]{Fren19DEGS},
\cite{PelSch22?CGST, FreMie21?DKRF}. However, we will give the main formal steps
to put the results into the perspective of Section \ref{su:Case2Product}. 

We first transform the problem as in Section with $\psi_\eps$ and $\psi_\eps$
from \eqref{eq:psieps.phieps}. With the notion from the previous subsection we
have $\ol\calE(U,w)=E(U)+ e(w)$ with 
\[
 E(U)=\int_{\Omega_\slow} \LB\big( \ol A(x)u(x)\big) \frac1{\ol A(x)}\, \dd x
 \quad \text{and} \quad 
 e(w)=\int_{\Omega_\fast} \LB\big( \ol A(x)w(x)\big) \frac1{\ol A(x)}\,\dd x
\]
(see \eqref{eq:EBz.R*cosh} for $\LB$) and the rescaled dual dissipation potential $\ol\calR^*(U,w;\Xi,\zeta) =
\calR_\slow^*(U,\Xi) + \calR_\fast^*(w,\zeta)+ \bm\delta_{\{0\}}\big(
\PPC_\slow \Xi{-} \PPC_\fast \zeta\big)$ with 
\[
\calR_\slow^*(U,\Xi)
 = \int_{\Omega_\slow}\!\!\!\! \frac{\ol K(Y)}2 |\pl_x\Xi(x)|^2 U(x)\dd x 
 \  \text{ and } \ 
\calR_\fast^*(w,\zeta) 
 = \int_{\Omega_\fast}\!\!\!\! \frac{\ol K(Y)}2 |\pl_x\zeta(x)|^2 w(x) \dd x.
\]

The reduced dissipation potential $\calR_\red$ is now obtained by the
saddle-point reduction, namely
\[
\frakB_{E,\calR_\red}(U,\Xi) = \sup_{w:\, P_\fast w=P_\slow U} \Big(
\inf_{\zeta:\, \PPC_\fast \zeta  =\PPC_\slow \Xi}  
\frakB_{e,\calR_\fast}(w,\zeta)\Big),
\]
where $\frakB_{e,\calR_\fast}(w,\zeta) :X_\fast\ti X_\fast^*\to \R$ takes the explicit form
(using $\rmD e(w) = \log(\ol Aw)$)
\begin{align*}
\frakB_{e,\calR_\fast}(w,\zeta) &= \calR^*_\fast(w,\zeta) -
\calR^*_\fast \big(w,-\rmD e(w)\big) 
\\
&= \int_{-1}^1 \frac12\Big(\ol K\,w\,\big|\pl_x\zeta|^2 - \frac{\ol K\,w}{(\ol
  A\,w)^2}\, \big|\pl_x(\ol A\,w)\big|^2     \Big) \dd x .
\end{align*}
It is surprising that the sup-inf of $\frakB_{e,\calR_\fast} $ under given boundary
conditions can be evaluated explicitly, see \cite[App.\,A]{LMPR17MOGG}, \cite[Sec.\,4]{Fren19DEGS}, and \cite[Sec.\,1.3]{PelSch22?CGST}. 

Here we provide a new and much shorter way of obtaining the desired result.

\begin{theorem}[Membrane reduction]
\label{th:MembrFormula}
Let $\ol K, \,\ol A \in \rmL^\infty([-1,1])$ be given and bounded from below by
a positive constant. Then 
\[
\calJ(w_-,w_+;\zeta_-,\zeta_+) : =\sup_{w(1)=w_+ \atop w(-1)=w_-} \ 
\inf_{\zeta(1)=\zeta_+ \atop \zeta(-1)=\zeta_-}  \frakB_{e,\calR_\fast}(w,\zeta)
\]
has the explicit form 
\begin{equation}
  \label{eq:MebranceExplForm}
\calJ(w_-,w_+;\zeta_-,\zeta_+)= K_\eff \sqrt{a_-w_-\,a_+w_+} \,
\sfC^*\big(\zeta_+{-}\zeta_- \big)  - K_\eff \, 2\big(\sqrt{a_+w_+}
- \sqrt{a_- w_-}\big)^2 
\end{equation}
where $K_\eff = \Big(\text{\small$\ds\int$}_{\hspace*{-0.41em}-1}^1 \ol
A(x)/\ol K(x) \,\dd x  \Big)^{-1} $, $a_+=\ol A(1^-)$, $a_-=\ol
A(-1^+)$,  and $\sfC^*$ is defined in \eqref{eq:EBz.R*cosh}. 
\end{theorem}
\begin{proof} Clearly, $ \frakB_{e,\calR_\fast}$ is strictly concave-convex and thus has at
  most one saddle point which is also the only critical point. Hence solving
  $\rmD \frakB_{e,\calR_\fast} = 0 $ gives the solution. 

However, it is advantageous to do a transformation first. We set
\[
w= v/\ol A, \quad  \zeta = \log(v/\eta^2) , \quad \text{and }\ 
\calI(v,\eta)= \frakB_{e,\calR_\fast} \big( v/\ol A , \log(v/\eta^2) \big).
\]
An elementary calculation shows that $\calI$ has a much simpler form, namely
\[
\calI(v,\eta) = -2 \int_{-1}^1 \ol\kappa \,\eta'\, \big(\frac v\eta\big)' \dd
x, \quad \text{where } \ol\kappa(x)=\ol K(x)/\ol A(x). 
\]
It will be particularly useful, that $\calI$ is linear in $v$. 

If $(w_*,\zeta_*)$ is a critical point for $\frakB_{e,\calR_\fast}$, then the
transformed point $(v_*,\zeta_*)$ is a critical point for $\calI$, and vice
versa. Hence, we have to determine the critical points of $\calI$ and observe
that
\[
\rmD_v \calI(v,\eta) = - \frac2\eta\:\big( \ol\kappa\, \eta'\big)'. 
\]
As the prefactor $2/\eta$ is irrelevant, we see that $\eta_*$ is uniquely
determined by its boundary values $\eta_-$ and $\eta_+$. In particular, we know
that $\ol\kappa \eta'_*$ must be constant, namely 
\[
 \ol\kappa(x) \eta'_*(x) = K_\eff \big(\eta_+ {-} \eta_-\big) \quad\text{for all }
 x \in [-1,1].
\]
Because of $\rmD_v \calI(v,\eta_*)=0$, this is enough to evaluate
$\calI(v,\eta_*)$ explicitly by only knowing the boundary values $v_-$ and
$v_*$ of $v_*$:
\[
\calI(u_*, \eta_*) = -2 \int_{-1}^1 \underbrace{\ol\kappa
  \,\eta'}_{=\text{const.}} \big( \frac{v}\eta\big)' \dd x = 2\;\!K_\eff
\big(\eta_+ {-} \eta_-\big) \, \Big( \frac{v_+}{\eta_+} 
 - \frac{v_-}{\eta_-} \Big). 
\]
Inserting the boundary conditions $ v_\pm = a_\pm w_\pm$ and $\eta_\pm = \big(
a_\pm w_\pm\big)^{1/2} \,\ee^{-\zeta_\pm/2}$ gives 
\[
\frakB_{e,\calR_\fast} ( w_*, \zeta_*) = - 2 K_\eff \Big( a_+ w_+ - \sqrt{a_+w_+ a_-w_-}
\big( \ee^{(\zeta_+-\zeta_-)/2} + \ee^{(\zeta_- - \zeta_+)/2} \big) + a_-w_-
\Big) ,
\]
which yields the desired formula \eqref{eq:MebranceExplForm}. 
\end{proof}

Using the port conditions $P_\fast w=P_\slow U$ and
$\PPC_\fast \zeta= \PPC_\slow \Xi$, the above result leads to the desired
BER structure
\[
\calJ(U_-,U_+;\Xi_-,\Xi_+)= \sfR_Y^*(U_-,U_+;\Xi_-,\Xi_+) - \sfR_y^*
\big(U_-,U_+;\log(A_-U_-),\log(A_+U_+)\big) ,
\]
where $A_-=\ol A(-1^-)$ and $A_+ = \ol A(1^+)$.

In summary, we obtain the effective gradient system $(X_\slow, E, \calR_\eff)$
with 
\[
\calR_\eff^*(U,\Xi)=\calR^*_\slow(U,\Xi) + K_\eff \sqrt{A_-U(-1)\,A_+U(1)} \:
\sfC^*\big(\Xi(1){-}\Xi(-1) \big).
\]
We see that the effective contribution of the membrane is of cosh-type, and in
particular it is not quadratic. Moreover, $\sfR_y$ depends on $\ol A$ which is
information that stems from $\calE_\eps$, which was not present in
$\calR_\eps^*$. Of course, also the cosh-type function $\sfC^*$ is inherited
from $\calE_\eps$, namely from the Boltzmann function $\LB$. Observe that
$ \mu=\LB'(r)=\log r$ has the inversion $r = \ee^\mu$. Using this for the
forward and backward fluxes it is no longer surprising to obtain $\sfC^*$.

Another way of understanding this transition involves the theory of
Markov processes. The pure diffusion problem can be interpreted as a particle
performing a Brownian motion with state-dependent mobility. In the membrane
region the mobility is very low (of order $\eps$) but the width of the membrane
is $2\eps$ such that the process will leave the membrane region very fast to
one or the other side. In the limit $\eps\to 0$ this gives rise to a jump
process for the particle either being reflected or jumping through the
barrier. According to \cite{MiRePe16GORR, MPPR17NETP} the corresponding
dissipation potential for such jump processes is defined via $\sfC^*$. 

\subsection{Linear reaction-diffusion equation}
\label{su:Linear.M.GS}

Before going into the one-dimensional equation with membrane scaling, we note
that the general structure of reaction-diffusion systems with detailed balance
condition has the following gradient structure. On $X=\rmL^1(\Omega;\R^{i_*})$
we consider 
\begin{align*}
&\calE(\bfc) = \calH(\bfc|\bfc_*)=\int_\Omega \LB(c_i/c^*_i) c^*_i\dd x \quad \text{and }
\\
&
\calR^*(\bfc;\bfxi) = \int_\Omega \Big(\sum_{i=1}^{i_*} \frac{K_i c_i}2 |\nabla
\xi_i|^2 + \sum_{r=1}^{r_*} \mu_r \big(\bfc^{\bfalpha^r}\bfc^{\bfbeta^r} \big)^{1/2}
\sfC^*\big(\bfxi{\vdot}(\bfalpha^r{-}\bfbeta^r ) \big)\Big) \dd x ,
\end{align*}
where $K_i\geq 0$ is the diffusion constants of species $X_i$, while $\mu_r>0$
is the reaction coefficient of the $r$th reaction having stoichiometric vectors
$\bfalpha^r,\, \bfbeta^r\in \N^{i_*}_0$. The associated gradient-flow equation
is the following system of $i_*$ equations:
\[
\dot c_i = \DIV\Big(K_i \big(\nabla c_i - \frac{c_i}{c^*_i}\nabla c^*_i\big) \Big)
- \sum_{r=1}^{r_*}  \mu_r \Big( 
\big(\tfrac{\ds\bfc_*^{\bfbeta^r}}{\ds\bfc_*^{\bfalpha^r}}\big)^{1/2} \bfc^{\bfalpha^r} -  
\big(\tfrac{\ds\bfc_*^{\alpha^r}}{\ds\bfc_*^{\bfbeta^r}}\big)^{1/2} \bfc^{\bfbeta^r} 
\Big) \big(\alpha_i^r - \beta_i^r\big) .
\]

In the same spirit as in the previous section we study again a linear PDE, but
now it has diffusion and reaction with the background, i.e.\
$A \leftrightharpoons \emptyset$. Again we assume that the material parameters
$K_\eps$ for diffusion and $B_\eps$ for reaction scale suitably with $\eps$ in
the membrane region ${]{-}\eps,\eps[}$. With
$\Omega_\eps ={]{-}1{-}\eps, 1{+}\eps [}$, the gradient system is given via
$X=\rmL^1(\Omega_\eps)$,
\begin{align*}
&\calE_\eps(u) =\calH(u|1/A_\eps) = \int_{\Omega_\eps} \LB(A_\eps u)
\frac1{A_\eps} \dd y \quad \text{and}
\\
&\calR^*(u,\xi)= \int_{\Omega_\eps}\Big( \frac{K_\eps u }2 |\xi'|^2 + B_\eps
\sqrt{u} \,\sfC^*(\xi) \Big) \dd y.
\end{align*}
Using $\phi_\eps:\Omega_\eps \to \Omega:=[-2,2]$ and
$\psi_\eps=\phi_\eps^{-1}:\Omega \to \Omega_\eps$ from \eqref{eq:psieps.phieps}
we assume that $A_\eps$, $B_\eps$, and $K_\eps$ are given in the form 
\begin{equation}
  \label{eq:RDS.coeff.scal}
  A_\eps(y)=\ol A(\psi_\eps(y)), \quad B_\eps(y)=\phi'_\eps(y) \,\ol
  B(\phi_\eps(y)), \quad K_\eps(y) = \frac1{\phi_\eps(y)} \,\ol K(\phi_\eps(y))
\end{equation}
for given functions $\ol A,\ \ol B, \, \ol K \in \rmP\rmC^0 \big(
[-2,-1]{\ol\cup}[-1,1]{\ol\cup}[1,2])$. To make our theory work we assume that
$\ol A$ and $\ol K$ have a positive lower bound on $\Omega$, whereas for $\ol
B$ it is sufficient to have $\ol B(x)\geq 0$. 

Transforming the system to the domain $\Omega$ as in the previous subsection,
we obtain a slow-fast gradient system $(X_\slow\ti X_\fast ,\calE_\eps,
\calR_\eps)$ given by  
\begin{align*} 
&X_\slow=\rmL^1(\Omega_\slow), \ \  X_\fast =\rmL^1(\Omega_\fast),  \ \
\Omega_\slow[-2,-1]\cup[1,2], \ \  \Omega_\fast=[-1,1],
\\
&\calE_\eps(U,w) = E(u)+\eps e(w),  \ \  E(U)=\int_{\Omega_\slow}\!\!\!
\LB(\ol A U)\frac{\dd x}{\ol A}, \ \ 
e(w)=\int_{\Omega_\fast}\!\!\! \LB(\ol A w)\frac{\dd x}{\ol A} ,
\\
&\calR^*_\eps(U,w; \Xi,\xi)=\ol\calR^*\big(U,w;\Xi,\frac1\eps\xi\big) \ \text{
  with } \\
&\ol\calR^*(U,w;\Xi,\zeta)=\ol\calR^*_\slow(U,\Xi) +
\ol\calR^*_\fast(w,\zeta) + \bm\delta_{\{\bm0\}}(P^0_\slow\Xi{-} \PPC_\fast
\zeta),
\\
&\ol\calR^*_\slow(U,\Xi) = \int_{\Omega_\slow} \Big(\frac{\ol K}{2}|\Xi'|^2 U +
\ol B\,\sqrt{U}\,\sfC^*(\Xi) \Big) \dd x, \text{ and } 
\\
&\ol\calR^*_\fast(w,\zeta) = \int_{\Omega_\fast} \Big(\frac{\ol K}{2}|\zeta'|^2 w +
\ol B\,\sqrt{w}\,\sfC^*(\zeta) \Big) \dd x.
\end{align*}

As in the previous subsection we obtain the effective gradient structure
$(X_\slow, E, \calR_\eff)$ by solving the sup-inf problem for the B-function
$\frakB_{\ol E,\ol\R}$ in the form $\calR_\eff^*=\calR_\slow^* + \calR_\red^*$
with $\calR^*_\red= \sfR^*_Y(P_\slow\cdot , \PPC_\slow)$, where we obtain an
explicit formula for $\sfR_Y$. To formulate the following result we introduce
the two auxiliary functions $H_+,\,H_-:[-1,1]\to \R$ via 
\begin{equation}
  \label{eq:Hpm.ODE}
  \big( \frac{\ol K}{\ol A} \,H'_\pm \big) ' = \frac{\ol B}{\ol A^{1/2}}\, H_\pm \
\text{ in } {]{-}1,1[}, \quad H_\pm(\pm1)=1, \quad H_\pm(\mp1)=0.
\end{equation}
Simple ODE arguments show $H_\pm(x)\in [0,1]$, $H'_-(x)<0$, and $H'_+(x)>0$ for
all $x\in [-1,1]$. 

\begin{theorem}[Membrane with reaction and diffusion]
 \label{th:MembRDS}
For the fast gradient system $(X_\fast,e, \calR^*_\fast)$ the reduced
B-function $\scrB_\red$ has the BER structure $(e,\sfR^*_Y)$ with 
\begin{align}
  \label{eq:MembRDS.red}
  \sfR^*_Y(w_-,w_+;\zeta_-,\zeta_+)
&= M_\eff \sqrt{\ol A(-1)w_-\ol A(1)w_+}\,\sfC^*(\zeta_+{-}\zeta_-) \\
\nonumber
& \quad   + M_- \sqrt{\ol A(-1)w_-} \,\sfC^*(\zeta_-)
   + M_+ \sqrt{\ol A(1)w_+}\, \sfC^*(\zeta_+),
\end{align}
where $M_\eff = \ol K(1) |H'_-(1)|/\ol A(1)=\ol K(-1)H'_+(-1)/\ol A(-1)$
and $M_\pm = \int_{-1}^1 \ol B \,H_\pm/\ol A^{1/2} \dd x$. 

In the case of constant coefficients we have 
\[
M_\eff = \frac{\ol K}{\ol A} \,\frac{\sigma\cosh(2\sigma)}{\sinh(2\sigma)}
\quad \text{and} \quad M_+=M_- = \frac{\ol B}{\ol A^{1/2}} \, { 
\frac{\cosh(2\sigma){-} 1}{\sigma \sinh(2\sigma)} \quad \text{with
}\sigma^2=\frac{\ol A^{1/2}\ol B}{\ol K}. }
\]
\end{theorem}
\begin{proof}
As in Theorem \ref{th:MembrFormula} we do a transformation to characterize the 
unique saddle point $(w_*,\zeta_*)$. With $w=v/\ol A$ and
$\zeta = \log (\ol A w/\eta^2)$, the B-function $\frakB_{e,\calR_\fast}$ gives 
\begin{equation}
  \label{eq:calI.MembRDS}
\calI(v,\eta):= \frakB_{e,\calR_\fast} \big(v/\ol A, \log(v/ \eta^2) \big)
 = \int_{-1}^1 \Big( {-}2 \ol\kappa \eta'\big(\frac v\eta\big)' +
2 \ol\beta \frac{1{-}\eta}\eta (v{-} \eta)  \Big) \dd y ,
\end{equation}
where $\ol\kappa = \ol K/\ol A$ and $\ol\beta=\ol B/\ol A{}^{1/2}$. Here we
used the specific interaction of $\sfC^*$ and $\log=\LB'$, namely $\sfC^*(\log
\alpha)= \big(\alpha^{1/4}{-}\alpha^{-1/4}\big)^2$. Of course, the construction
is such that $\eta\equiv 1$ leads to $\calI(v,1)=0$. 

The surprising and helpful fact is that $\calI$ is affine in $v$ which allows
us to evaluate $ \frakB_{e,\calR_\fast} (w_*,\zeta_*)=\calI(v_*,\eta_*)$ at the
unique critical point. In particular, we have
\[
 0 = \rmD_v \calI(v,\eta) = \frac2\eta \Big( \big(\ol\kappa \eta'\big)' -
 \ol\beta\,\eta + \ol\beta\Big),
\] 
such that the critical point $(v_*,\eta_*)$ satisfies the linear ODE 
$- (\ol\kappa \eta')'+\ol\beta \eta= \ol \beta$. Hence, 
\begin{align*}
\frakB_{e,\calR_\fast} (w_*,\zeta_*)&= \calI(v_*,\eta_*) =  \int_{-1}^1
 \Big( {-}2 \ol\kappa \eta_*'\big(\frac{ v_*}{ \eta_*} \big)' +
2 \ol\beta \frac{1{-}\eta_*}{\eta_*} (v_*{-} \eta_*)  \Big) \dd y 
\\
&=\Big[{-}2 \ol\kappa \eta'_*\frac{v_*}{\eta_*} \Big]_{x=-1}^1 + \int_{-1}^1 \Big(
 \frac2{\eta_*}\big(\underbrace{ (\ol\kappa \eta'_* )' +\ol\beta(1{-} \eta_*)
 }_{=0} \big) + 2\underbrace{\ol\beta(\eta_*{-}1)}_{= (\ol\kappa \eta'_*)'} \Big) \dd x 
\\
& = 2\ol\kappa_-\eta'_*(-1) \big(\frac{v_-}{\eta_-} -1\big) + 
2\ol\kappa_+ \eta'_*(1)\big(1-\frac{ v_+}{\eta_+}\big),
\end{align*}
where $\ol\kappa_\pm=\ol\kappa(\pm1)$, $v_\pm=v_*(\pm1)$, and
$\eta_\pm=\eta_*(\pm1)$. 

Using the auxiliary functions $H_\pm$ we have $\eta_*=1+(\eta_-{-}1)H_- +
(\eta_+{-}1) H_+$ which gives $ \eta'_*(\pm 1)= (\eta_- {-}1) H'_-(\pm 1)
+(\eta_+{-}1) H'_+(\pm 1) $. Abbreviating $b_\pm:=\sqrt{v_\pm}$ and
$E_\pm:=\ee^{\zeta_\pm/2}$ and using $\eta_\pm = \sqrt{v_\pm}\ee^{-\zeta_\pm/2}
= b_\pm E_\pm^{-1}$ we obtain
\begin{align*}
\frakB_{e,\calR_\fast} (w_*,\zeta_*)
&= 2\ol\kappa_+ H'_+(1)\big( b_+(E_+{+} E_+^{-1}){-} b_+^2{-}1\big)
%% \\&\quad 
- 2\ol\kappa_-H'_-(-1)\big( b_-(E_-{+} E_-^{-1}){-}b_-^2{-}1 \big)
\\
&\quad +2\ol\kappa_- H'_+(-1)(b_+E_+^{-1}{-}1)(b_-E_-{-}1) 
%% \\ &\quad 
 -2\ol\kappa_+ H'_- (1)(b_-E_-^{-1}{-}1) (b_+E_+{-}1). 
\end{align*}
To simplify this expression, we use that the Wronski determinant 
$\ol\kappa H'_+H_-- \ol\kappa  H'_-  H_+$ is constant on $[-1,1]$, and we call
this constant $M_\eff>0$. Using the boundary conditions of $H_\pm$ we have 
$M_\eff= \ol\kappa_- H'_+(-1)= -\ol\kappa_+ H'_-(1)$. Moreover, integrating the
ODE \eqref{eq:Hpm.ODE} yields
\[
\pm \ol\kappa_\pm H'_\pm(\pm1) = \pm \ol\kappa_\mp H'_\pm(\mp1) +\int_{-1}^1
\ol\kappa H_\pm \dd x =M_\eff + M_\pm ,
\]
by exploiting our definition of $M_\pm$ in \eqref{eq:MembRDS.red}. 
With this we arrive at 
\begin{align*}
\frakB_{e,\calR_\fast} (w_*,\zeta_*)
&= 2 M_- \big( b_- (E_-{+} E_-^{-1}{-}2) - (b_-{-}1)^2\big)
%\\ &\quad  
+ 2M_+ \big( b_+ (E_+{+} E_+^{-1}{-}2) - (b_+{-}1)^2\big) 
\\
&\quad +2M_\eff \big( b_+b_-(E_+E_-^{-1}{+}E_+^{-1}E_-{-}2) -(b_+{-}b_-)^2\big).
\end{align*}
Inserting $E_\pm = \ee^{-\zeta_\pm/2}$ and $b_\pm= \sqrt{v_\pm}=\sqrt{a_\pm
  w_\pm}$ yields 
\begin{align*}
\frakB_{e,\calR_\fast} (w_*,\zeta_*)&=
 M_-\big( \sqrt{a_-w_-}\,\sfC^*(\zeta_-) - 2(\sqrt{a_-w_-}{-}1)^2\big)
\\ 
 &\quad  
  +M_+\big( \sqrt{a_+w_+}\,\sfC^*(\zeta_+) - 2(\sqrt{a_+w_+}{-}1)^2\big) 
\\ &\quad + M_\eff\,\big( \sqrt{a_+w_+a_-w_-}\,\sfC^*(\zeta_+{-}\zeta_-) -
2(\sqrt{a_+w_+}{-}\sqrt{a_-w_-})^2 \big)
\\
&= \sfR_Y^*(w_-,w_+;\zeta_-,\zeta_+) - \sfR_Y^*\big(w_-,w_+;
\log(a_-w_-),\log(a_+w_+)\big), 
\end{align*}
which is the desired general formula \eqref{eq:MembRDS.red}.

The special formula for constant coefficients follows by setting $\sigma^2= \ol
A^{1/2} \ol B/ \ol K$ and observing $H_\pm(x) =\sinh(\sigma \pm \sigma
x)/\sinh(2\sigma)$. 
\end{proof}

\appendix
\section{Classical existence theory for saddle points}
\label{se:Saddles}

We recollect the basic result from saddle point theory as contained in
\cite[Cha.\,VI]{EkeTem74ACPV} (La dualit\'e par les minimax).

We consider a functional $\scrL:\bfU\ti \bfV \to \ol\R=[-\infty,\infty]$, 
called ``Lagrangian'' in \cite[Cha.\,VI]{EkeTem74ACPV}, where we now want
to minimize with respect to $x\in \bfU$ and maximize with respect to
$y \in \bfV$. This means that for applying the theory below to the 
B-functions used above we have to set $ \bfU=X$, $\bfV=X^*$, and
$\scrL(u,\xi)=- \scrB(u,\xi)$. Now, a point $(x_*,y_*)$ is called a \emph{saddle
  point} of $\scrL$ if
\[
\forall \,x\in \bfU,\ y\in \bfV:\quad \scrL(x_*,y) \leq \scrL(x_*,y_*) \leq
\scrL(x,y_*). 
\]
Thus, we minimize with respect to $x\in \bfU$, and we maximize with
respect to $y$. 

The aim is to find a saddle point from general principles. For this one looks
at $\sup_{y\in \bfV} \inf_{x\in \bfU} \scrL(x,y)$ and
$\inf_{x\in \bfU} \sup_{y\in \bfV} \scrL(x,y)$. We obviously always have a one-sided
estimate, and the major question in constructing saddle points is when we have
equality.

\begin{lemma}[Simple facts on saddles points] 
\label{le:SimpleFactsSP}
\begin{align}
\text{(a)}\qquad &  \label{eq:SimpleEst}
 \SI_\scrL:= \sup_{y\in \bfV} \inf_{x\in \bfU} \scrL(x,y) \ \leq \ 
   \inf_{x\in \bfU} \sup_{y\in \bfV} \scrL(x,y) := \IS_\scrL
\\
\text{(b)}\qquad & \label{eq:Sadd.equal}
\text{saddle point $(x_*,y_*)$ exists} \ \ \Longrightarrow \ \ 
\SI_\scrL=\IS_\scrL. 
\end{align}
In the latter case, we have $\scrL(x_*,y_*)=\SI_\scrL=\IS_\scrL$. 
\end{lemma}
\begin{proof}
To show (a), we start from $\scrL(x,y)\leq \sup_{\ol y} \scrL(x,\ol
y)$. Taking the infimum over $x$ we obtain $\inf_{x} \scrL(x,y) \leq
\IS_\scrL$. Now taking the supremum over $y$ in the left-hand side
leads to the desired estimate $\SI_\scrL \leq \IS_\scrL$.

To show (b) simply note that the saddle-point property implies 
\[
\inf_{x\in \bfU} \scrL(x,y_*)=\scrL(x_*,y_*) = \sup_{y\in \bfV}
\scrL(x_*,y).
\]
Thus, we find $\SI_\scrL\geq \scrL(x_*,y_*) \geq \IS_\scrL$. With (a)
this implies the desired equality. 
\end{proof}

The quantity $ \delta_\scrL= \IS_\scrL- \SI_\scrL \geq 0$ is called the \emph{duality
gap}. The function $\scrL(x,y)=\tanh(x{-}y)$ on $\R\ti \R$ shows that
$\delta_\scrL$ can be positive. Indeed, $\SI_{\tanh}=-1$ and $\IS_{\tanh}=+1$ such
that $\delta_{\tanh} = 2$. 

The opposite implication in \eqref{eq:Sadd.equal} is not
valid. To see this consider $\bfU=\bfV=\R$ and
$\scrL(x,y)=\ee^x-\ee^{-y}$. Clearly, $\inf_x \scrL(x,y)= - \ee^{-y}$
and hence, $\SI_\scrL=0$ and similarly $\IS_\scrL=0$. However, no
saddle-point exists. Even in cases where no saddle-point exists it is an interesting
question under what conditions the duality gap is $0$, see e.g.\
\cite[Ch.\,III,\,Prop.\,2.3]{EkeTem74ACPV}. 

If two saddle points $(x_j,y_j)$ with $j=1,2$ exist, we have
\[
\scrL(x_1,y_2) \leq  \scrL(x_1,y_1) \leq  
\scrL(x_2,y_1) \leq  \scrL(x_2,y_2) \leq  \scrL(x_1,y_2) ,
\]
which means that all four points have the same value. If each
$\scrL(\cdot, y_j)$ is convex and each $\scrL(x_j,\cdot) $ concave,
then we conclude $\scrL(x,y)=\scrL(x_1,y_1)$ for all
$x=(1{-}s)x_1+sx_2$ and $y=(1{-}r)y_1+ry_2$ with arbitrary $r,s\in
[0,1]$. 
\medskip

A standard existence result for saddle points can be found in
\cite[Ch.\,VI,\,Prop.\,2.1]{EkeTem74ACPV}. We provide a variant that is
adjusted to our purposes.

\begin{proposition}[Existence of saddle points]
\label{pr:ExiSaddlePoint}
Consider reflexive Banach spaces $\bfU$ and $\bfV$ and assume that the
following conditions hold:
\begin{subequations}
\begin{align}
&\label{eq:Cond.a}
 \forall\, y\in \bfV: \quad x\mapsto \scrL(x,y) \text{ is convex and lsc},\\
&\label{eq:Cond.b}
 \forall\, x\in \bfU: \quad y\mapsto -\scrL(x,y) \text{ is convex and lsc},\\
&\label{eq:Cond.c}
\exists\,y_0\in \bfV:\quad \scrL(\,\cdot\,, y_0)  \text{ is coercive},\\
&\label{eq:Cond.d}
\exists\,x_0\in \bfU:\quad \!{-}\scrL(x_0,\,\cdot\,)  \text{ is coercive}.
\end{align}
\end{subequations}
Then, a saddle point $(x_*,y_*)$ exists and 
\[
\scrL(x_*,y_*)= \min_{x\in \bfU} \sup_{y\in \bfV} \scrL(x,y) = 
\max_{y\in \bfV} \inf_{x\in \bfU}  \scrL(x,y). 
\]
If moreover, in \eqref{eq:Cond.a} and \eqref{eq:Cond.b} we have strict
convexity, then the saddle point is unique. 
\end{proposition} 
\begin{proof} 
\STEP{Step 1: Saddle points on balls using strict convexity.} 
We additionally impose that
\begin{equation}
  \label{eq:calL.strcvx}
  \forall\, y\in \bfV:\quad \scrL(\,\cdot\,,y):\bfU\to \ol\R \text{ is strictly convex}.
\end{equation}
For $R\geq R_0:=
\max\{\|x_0\|_\bfU,\| y_0\|_\bfV\}$ we consider the closed and convex balls
$\bfU_R=\set{x\in \bfU}{\| x\|_\bfU\leq R}$ and similarly $\bfV_R$. 

For all $R$ we obtain a saddle point $(x_R,y_R)$ as follows.  For all
$y \in \bfV_R$ the direct method of the calculus of variations provides a
minimizer $x=\wh x_R(y) \in \bfU_R$ for $\scrL(\cdot,y)|_{\bfU_R} $, i.e.\
$\scrL(\wh x_R(y),y)=\min_{x\in \bfU_R} \scrL(x,y)=:\lambda_R(y)$. By the strict
convexity in \eqref{eq:calL.strcvx} $\wh x_R(y)$ is uniquely determined.

We first observe that $-\lambda_R:\bfV_R \to \ol\R$ is convex and lsc, as it is
the supremum of the convex and lsc functions $-\scrL(x,\cdot)$. Moreover, by
\eqref{eq:Cond.d} the function $-\lambda_R$  is
bounded from below by the proper, lsc, convex function $-\scrL(x_0,\cdot)$. 
Hence, $\lambda_R$ attains its maximum in a point $y^R \in \bfV_R$.  

Our aim is now to show that $(\wh x_R(y^R),y^R)$ is a saddle point of $\scrL$
on $\bfU_R\ti \bfV_R$. For this we choose arbitrary $y \in \bfV_R$ and $\theta\in [0,1]$
and set $x_\theta(y):=\wh x_R((1{-}\theta)y^R{+}\theta y)$ and obtain
\begin{align*}
\lambda_R(y^R) &\geq \lambda_R\big((1{-}\theta)y^R{+}\theta y \big) =
\scrL\big(x_\theta(y), (1{-}\theta)y^R{+}\theta y  \big) \\
&\overset{-\scrL(x_\theta(y),\cdot)\text{ cvx}}\geq (1{-}\theta)\scrL(x_\theta,
y^R) + \theta \scrL(x_\theta(y),y) \geq (1{-}\theta)\lambda_R(y^R) + \theta
\scrL(x_\theta(y),y). 
\end{align*}
In particular, for $\theta \in{]0,1]}$ and all $y\in \bfV_R$ we conclude 
\begin{equation}
  \label{eq:EstimSP.A5}
\lambda_R(y^R) =\scrL(\wh x_R(y^R))\geq \scrL(x_\theta(y),y) \quad \text{ for
  all }y \in \bfV_R.
\end{equation}  

Choosing $\theta=1/k$ for $k\in \N$, we obtain $x_k:= x_{1/k}(y)\in \bfV_R$ and
may select a weakly convergent subsequence (not relabeled) with $x_k \weak
x^R$. We claim that $x^R=\wh x_R(y^R)$ and hence is independent of
$y$. Indeed, for our fixed $y
\in \bfV_R$ and arbitrary $\wt x\in \bfU_R$ we have 
\begin{align*}
\scrL(x^R,y^R) &
\overset{\text{\eqref{eq:Cond.a},lsc}}\leq \liminf_{k\to \infty}\scrL(x_k,y^R)  
\\
& \overset{\text{\eqref{eq:Cond.b},cvx}}\leq 
 \limsup_{k\to \infty} \frac1{1{-}\frac1k}\Big(
 \scrL(x_k,(1{-}\tfrac1k)y^R{+}\tfrac1ky) - \frac1k\scrL(x_k,y)\Big)
\\
&\overset{\text{def.\,$\lambda_R,\wh x_R$}}\leq \limsup_{k\to \infty} \Big(
\frac{k}{k{-}1} 
 \scrL(\wt x,(1{-}\tfrac1k)y^R{+}\tfrac1ky) - \frac1{k{-}1} \lambda_R(y)\Big)
\\
&\overset{\lambda_R(y)<\infty}\leq \limsup_{k\to \infty}
\scrL(\wt x,(1{-}\tfrac1k)y^R{+}\tfrac1ky) 
\
\overset{\text{\eqref{eq:Cond.b},lsc}}\leq \ \scrL(\wt x, y^R),
\end{align*}
where we used $(1{-}\tfrac1k)y^R{+}\tfrac1ky \to y^R$ in the last step. Since
$\wt x\in \bfU_R$ was arbitrary we obtain $\lambda_R(y^R)\leq \scrL(x^R,y^R)\leq
\min_{\wt x \in \bfU_R} \scrL(\wt x,y^R)=\lambda_R(y^R)$. Hence, $\bfU^R$ is a
minimizer of $\scrL(\cdot,y^R)$ and hence coincides with $\wh x_R(y^R)$ because
of the strict convexity \eqref{eq:calL.strcvx}.    

Because of the uniqueness of the limit we conclude that for all $y\in \bfV_R$ we
have $x_\theta(y) \weak x^R=\wh x_R(y^R)$ for $\theta\to 0^+$. Thus, taking the
limit  $\theta\to 0^+$ in \eqref{eq:EstimSP.A5} and exploiting the lsc from
\eqref{eq:Cond.a} we obtain 
\[
\forall\, y\in \bfV_R\ \forall \, \wt x \in \bfU_R:\quad \scrL(x^R,y) \leq
\lambda_R(y^R)=\scrL(x^R,y^R) \leq \scrL(\wt x , y^R).
\]
This shows that $(x^R,y^R)$ is a saddle point for $\scrL$ restricted to $\bfU_R\ti
\bfV_R$. 

\STEP{Step 2: Saddle points on balls without strict convexity.} If we only have
convexity we consider 
\[
\scrL_\eps(x,y)=\scrL(x,y) + \eps \|x\|^2 \quad \text{with } \eps >0,
\]
where we can choose a strictly convex norm $\|\cdot\|$ on the reflexive space
$\bfU$. By Step 1 we obtain a saddle point $(x^R_\eps,y^R_\eps)\in \bfU_R\ti
\bfV_R$. Hence, we have 
\begin{equation}
  \label{eq:eps.saddle}
  \forall\, y\in \bfV_R\ \forall \, x \in \bfU_R:\quad \scrL(x^R_\eps,y)+
\eps\|x^R_\eps\|^2 \leq \scrL(x^R_\eps,y^R_\eps) + \eps\|x^R_\eps\|^2 
\leq \scrL(x , y^R_\eps)+ \eps \| x\|^2.
\end{equation}

We may choose a subsequence (not relabeled) with $(x^R_\eps,y^R_\eps) \weak
(\ol x^R, \ol y^R)$ in $\bfU\ti \bfV$. Dropping the middle term in
\eqref{eq:eps.saddle} we can pass to the limit using the lsc in \eqref{eq:Cond.a} and
\eqref{eq:Cond.b} and arrive at 
\[
  \forall\, y\in \bfV_R\ \forall \, x \in \bfU_R:\quad \scrL(\ol x^R,y) \leq \scrL(x , \ol y^R).
\] 
Thus, $(\ol x^R,\ol y^R)$ is indeed a saddle point for $\scrL$ restricted to
restricted to $\bfU_R\ti \bfV_R$. 

\STEP{Step 3: Unbounded case.} We now consider the limit $R\to \infty$. 
Using the coercivities \eqref{eq:Cond.c} and \eqref{eq:Cond.d}. For $R\geq R_0$
the saddle points $(x^R,y^R)$ from Step 2 satisfy 
\begin{equation}
  \label{eq:Saddle3}
  \scrL(x_R,y_0)\leq \scrL(x^R,y^R) \leq \scrL(x_0,y^R).
\end{equation}
Since $\scrL(\cdot,y_0)$ and $-\scrL(x_0, \cdot)$ are lsc and coercive (cf.\
\eqref{eq:Cond.c} and \eqref{eq:Cond.d}), they are bounded from below:
\[
\exists\, M>0\ \forall\, x \in \bfU\ \forall\,y\in \bfV:\quad  
\scrL(x,y_0)\geq -M \ \text{ and } \ \scrL(x_0,y)\leq M.
\]  
Combining this with \eqref{eq:Saddle3}, we have 
\[
\forall\, R\geq R_0:\quad  \text{(i)  }\scrL(x^R,y_0)  \geq -M \ \text{ and } \
\text{(ii) } \scrL(x_0,y^R)\leq M. 
\]
With \eqref{eq:Saddle3} we obtain $|\scrL(x^R,y^R)|\leq M$. 
Using the coercivity \eqref{eq:Cond.c} and (ii) we find $\|y^R\|\leq C_\bfV$, and
similarly \eqref{eq:Cond.d} and (i) give $\|x^R\|\leq C_\bfU$. Thus, using the
reflexivity of $\bfU$ and $\bfV$ we find a subsequence $(x^R,y^R)$ (not relabeled) such that 
\[
\scrL(x^R,y^R)\to \lambda_*,\quad x^R\weak x_* \text{ in }\bfU, \quad 
y^R\weak y_* \text{ in }\bfV.
\] 

For arbitrary $x \in \bfU$ we choose $R>\max\{R_0,\|x\| \}$ and obtain 
$\scrL(x,y^R) \geq \scrL(x^R,y^R) $. Taking the limit $R\to \infty$ (along the
subsequence) and using the lsc of $-\scrL(x,\cdot)$ we arrive at
\[
\scrL(x,y_*) \geq \limsup_{R\to \infty} \scrL(x,y^R)  \geq   \limsup_{R\to
  \infty}\scrL(x^R,y^R) =\lambda_*, 
\]
where $x\in \bfU$ was arbitrary. Similarly, we obtain $\scrL(x_*,y)\leq \lambda_*$
which gives the desired saddle-point property for $(x_*,y_*)\in \bfU\ti \bfV$:
\[
\forall\, x\in \bfU,\ y\in \bfV:\quad \lambda(x_*,y)\leq
\lambda_*=\lambda(x_*,y_*)\leq \scrL(x,y_*).
\]

\STEP{Step 4: Uniqueness under strict convexity.} This was shown already in Step
1.  

This completes the proof of Proposition \ref{pr:ExiSaddlePoint}. 
\end{proof}

\paragraph*{Acknowledgments.} This research was partially funded by Deutsche
Forschungsgemeinschaft (DFG) through the Collaborative Research Center SFB 1114
``{\itshape Scaling Cascades in Complex Systems}'' (Project no.\ 235221301)
within the Subproject C05 ``Effective models for materials and interfaces with
multiple scales''.

{\footnotesize
%\bibliographystyle{my_alpha}
%\bibliography{alex_pub,bib_alex}

\newcommand{\etalchar}[1]{$^{#1}$}
\def\cprime{$'$}

}

\end{document}